\documentclass[a4paper]{article}
\usepackage[T1]{fontenc}
\usepackage{amsmath}
\usepackage{amsthm}
\usepackage{amssymb}
\usepackage{graphicx}
\usepackage{authblk}
\usepackage{hyperref}
\usepackage{multirow}
\usepackage{subfigure}
\newcommand{\argmax}{\operatorname{argmax}}
\usepackage{cite}
\usepackage{color} 
\usepackage[dvipsnames]{xcolor}
\usepackage{comment}
\usepackage{enumitem}
\newcommand{\coeff}{\operatorname{coeff}}
\newcommand{\degree}{\operatorname{deg}}
\newcommand{\Conv}{\mathbf{trop}}
\newtheorem{theorem}{Theorem}[section]

\newtheorem{lemma}[theorem]{Lemma}

\newtheorem{problem}{Problem}
\newtheorem{definition}{Definition}
\newtheorem{remark}{Remark}

\newcommand{\e}{{\mathrm e}}
\usepackage{fullpage}
\newcommand{\eqlab}[1]{\label{eq:#1}}
\renewcommand{\eqref}[1]{(\ref{eq:#1})}

\newcommand{\Wert}{{\vert\kern-0.25ex\vert\kern-0.25ex\vert}}

\newcommand{\m}{{\operatorname{m}_\epsilon(\alpha,\gamma)}}

\newcommand{\md}{{\mathbf{d}}}
\newcommand{\figref}[1]{Fig.~\ref{fig:#1}}
\newcommand{\figlab}[1]{\label{fig:#1}}

\newcommand{\lemmaref}[1]{Lemma~\ref{lemma:#1}}
\newcommand{\lemmalab}[1]{\label{lemma:#1}}
\newcommand{\remref}[1]{Remark~\ref{remark:#1}}
\newcommand{\remlab}[1]{\label{remark:#1}}
\newcommand{\thmref}[1]{Theorem~\ref{theorem:#1}}
\newcommand{\thmlab}[1]{\label{theorem:#1}}

\newcommand{\defnlab}[1]{\label{defn:#1}}
\newcommand{\defnref}[1]{Definition~\ref{defn:#1}}

\newcommand{\secref}[1]{Section~\ref{sec:#1}}
\newcommand{\seclab}[1]{\label{sec:#1}}

\newcommand\rsp[1]{{\color{black}{#1}}}

\author[1]{Kristian Uldall Kristiansen}
\author[2]{Peter Szmolyan}
\affil[1]{Technical University of Denmark, Lyngby, 2800 Kgs. Lyngby, Denmark}
\affil[1]{{\tt krkri@dtu.dk}}
\affil[2]{TU Vienna, Austria}
\affil[2]{{\tt peter.szmolyan@asc.tuwien.ac.at}}

\title{Tropicalization and slow manifolds: A case study of the Michaelis--Menten reaction}

\begin{document}
\maketitle
\date{}

\begin{abstract}\noindent
We present a case study of a recently developed approach to multiple-scale systems based on exponential scaling and tropicalization (Kristiansen and Sarantaris, \emph{On a tropicalization of planar polynomial ODEs}). Parameters and variables are scaled exponentially with respect to a small parameter $\epsilon$, so that orders of magnitude become a genuine non-smooth singular limit $\epsilon\to 0$. Tropical geometry, via subdivisions of the Newton polygon, classifies the asymptotic regimes and locates the candidates for slow manifolds on all scales of the variables simultaneously. Here we combine this framework with blow-up and geometric singular perturbation theory. {In particular,} the non-smooth singular limit is resolved by spherical and cylindrical blow-ups, and degeneracies caused by exponentially small terms with flat contact are treated by adding the exponential quantities as new variables {(following Kristiansen, \emph{Blowup for flat slow manifolds})}.
As the case study, we take the Michaelis--Menten reaction. The tropical classification yields three asymptotic regimes of parameters. In two of these, we find a sequence of blowups and desingularizations where the system has a {global} normally hyperbolic critical manifold on the blown up space. In the remaining case, we obtain a similar result, but only for a subcase. In particular, there is an asymptotic regime of parameters where there is no global slow manifold of the Michaelis--Menten reaction. 
\end{abstract}
\textit{Keywords:} Tropicalization, GSPT, blow-up, reaction-kinetics, exponential asymptotics
\newline
\textit{2020 Mathematics Subject Classification:} 34C23, 34C45, 37G05, 37G10


\section{Introduction}


Polynomial systems of ordinary differential equations (ODEs) occur in many different areas of science, perhaps most notably as models of chemical reaction networks \cite{feinberg2019a}. Frequently, dependent variables and parameters vary on vastly different scales, leading to multiscale behaviour of solutions. In suitably scaled forms of the equations, dimension reduction methods can be used to obtain qualitative insight. A prime example is dimension reduction by Quasi-Steady-State-Assumption (QSSA), which is based on the assumption that certain concentrations (involved in fast reactions) are equilibrated for an extended (slow) time \cite{segel1989a}. One way of making QSSA precise is through the concept of slow manifolds from Geometric Singular Perturbation Theory (GSPT). These objects are locally invariant manifolds capturing the essential dynamics after initial transients \cite{fen3}. GSPT, in combination with the technical tool of blowup, is by now a well-developed toolbox for systems depending singularly on a distinguished small parameter \cite{kuehn2015,wechselberger2020a}. In different regions of phase space and parameter space, different types of slow manifolds can exist in separate scales. However, finding the relevant scales is often challenging \cite{noel2014a}.

Recently, concepts from tropical geometry have been used to find the governing scales in a systematic and algorithmic way, see \cite{desoeuvres2022a,noel2014a,soliman2014a}. The basic idea is to scale all parameters and variables according to certain weights and then determine the weights such that dominant terms of opposite signs cancel. This gives candidates for dimension reduction (in the spirit of slow manifolds). In \cite{kristiansen2026a} this approach was formalized for planar systems by writing the parameters and variables in suitable logarithmic coordinates. In a singular limit, one obtains a piecewise smooth system with tropical curves as discontinuity sets. Sliding along suitable tropical curves is a proxy for motion on slow manifolds, see also \cite{kristiansen2026a,portisch2021novel}.

In this paper, we relate (for the first time) the dynamics on the tropical level (which is singular and discontinous) back to the polynomial system by applying methods from GSPT. 
We do this in the context of the Michaelis--Menten reaction model:
\begin{equation*}
E + S \;\underset{k_{-1}}{\overset{k_1}{\rightleftharpoons}}\; C \;\overset{k_2}{\longrightarrow}\; E + P.
\end{equation*}
described by the following equations
\begin{equation}\eqlab{mm-dim}
\begin{aligned}
\frac{ds}{dt} &= -k_1(e_0 - c)\,s + k_{-1}c,\\
\frac{dc}{dt} &= k_1(e_0 - c)\,s - (k_{-1} + k_2)\,c,
\end{aligned}
\end{equation}
for the substrate $s$ and complex $c$.  Here $e_0$ denotes the total enzyme concentration while $k_{-1}$, $k_1$ and $k_2$ are the rate constants.
We will work with a dimensionless form of the equations (following \cite[Eq. 9]{lapuz}) based upon the following scalings, where $s_0$ represents the initial substrate concentration:
\begin{equation}\eqlab{tueqn}
\tilde u = \frac{s}{s_0},\quad \tilde v = \frac{c}{e_0},\quad \tau = k_1 e_0 t,\qquad
\tilde\alpha = \frac{k_{-1}}{k_1 s_0},\quad \tilde\beta = \frac{e_0}{s_0},\quad \tilde\gamma = \frac{k_2}{k_1 s_0}.
\end{equation}
Here $\tilde u$ and $\tilde v$ denote the nondimensionalized substrate and complex while $\tilde \alpha$, $\tilde \beta$ and $\tilde \gamma$ are nondimensionalized reaction rates. This brings \eqref{mm-dim} with $' = d/d\tau$ into the following form:
\begin{equation}\eqlab{mm}
\begin{aligned}
 \tilde{u}' &=- \tilde{u} + \tilde{\alpha} \tilde{v}+ \tilde{u}\tilde{v},\\
 \tilde \beta\tilde{v}'&=\tilde{u}-(\tilde{\alpha}+\tilde{\gamma})\tilde{v}- \tilde{u}\tilde{v}.
\end{aligned}
\end{equation}
Up to the rescaling of time by $\tilde \beta$, this is the form used in \cite{lapuz}
%

We recognize that there are different ways of non-dimensionalizing the Michaelis--Menten reaction model, but {in this work} we work ``globally'' in all of $(\tilde u,\tilde v)\in (0,\infty)^2$ so these differences will not be important.

The study of \eqref{mm} as a multiscale problem and the associated celebrated quasi-steady-state approximation \cite{michaelis1913a,briggs1925a} has been a guiding example for the theory of singular perturbations. Singular perturbation analyses start with \cite{heineken1967a}, and the identification of the relevant small parameters by scaling and time scale estimates with \cite{segel1988a,segel1989a}; see \cite{borghans1996a,schnell2000a,tzafriri2003a} for the total and the reverse quasi-steady-state approximation. A systematic and coordinate-independent approach, based on Tikhonov--Fenichel parameter values, was developed by Walcher and coworkers \cite{noethen2011a,goeke2014a,goeke2015a,lax2020a}.
It has been applied to \eqref{mm} in \cite{eilertsen2020a}. A reduced model based on computational singular perturbation, together with local validity criteria, is given in \cite{patsatzis2019a,patsatzis2023a}, and a catalogue of reductions by geometric singular perturbation theory in \cite{lapuz}. From the rigorous point of view, quasi-steady-state approximations can often be justified through GSPT, as critical manifolds perturbing to nearby slow manifolds. For technical reasons, these manifolds are often local {in phase space}. A comparison of our findings with existing work is given in the discussion, see \secref{conclude}.

It is a simple fact that the equilibrium $(\tilde u,\tilde v)=(0,0)$ of \eqref{mm} attracts all points in the first quadrant.  Moreover, it is well-known (\cite{davis1999a,calder2008a}) that there is unique invariant manifold from infinity for all parameters $(\tilde \alpha,\tilde \beta,\tilde \gamma)$. It has the following asymptotics:
\begin{align*}
 \tilde v = 1+\mathcal O(\tilde u^{-1}),\quad \tilde u\rightarrow \infty.
\end{align*}
This also follows directly from a compactification of $\tilde u\ge  k$, $k>0$, by $\tilde z=\tilde u^{-1}$ so that
\begin{equation}\eqlab{zv}
\begin{aligned}
 \tilde{z}' &=\tilde z^2 \left(1-\tilde \alpha \tilde z \tilde v  -\tilde v\right),\\
 \tilde \beta \tilde{v}'&=1-(\tilde{\alpha}+\tilde{\gamma})\tilde z \tilde{v}- \tilde{v}.
\end{aligned}
\end{equation}
after multiplication of the right hand side by $\tilde z$. The invariant manifold from infinity then corresponds to the unique center manifold of the partially hyperbolic point $(\tilde z,\tilde v)=(0,1)$. (The reference \cite{eilertsen2021a} identifies a similar trajectory for \eqref{mm} with a constant inflow in the $\tilde u$-equation.)
%
Frequently, the invariant manifold from infinity is referred to as a slow manifold without any rigorous justification. Classical results of GSPT provide slow manifolds on compact subsets of phase space for one fixed scaling of the variables, i.e. locally. 

We use the following problem as a test case of our approach on combining tropical methods with GSPT and blowup:
\begin{problem}\label{problem}
\textcolor{black}{In what regimes of the $(\tilde \alpha,\tilde \beta,\tilde \gamma)$-parameter space can the unique manifold from infinity be understood as a global slow manifold obtained by perturbing a well-defined critical manifold?}
\end{problem}

In particular, we are concerned with a  characterization of all generic regimes in the $(\tilde \alpha,\tilde \beta,\tilde \gamma)$-parameter space, going beyond a discrete set of scalings by powers of a small parameter $\epsilon$, such as $\tilde d=\mathcal O(\varepsilon),\mathcal O(1)$, or $\mathcal O(\varepsilon^{-1})$ for all $d\in \{\alpha,\beta,\gamma\}$ as $\varepsilon\rightarrow 0$, see \cite{lapuz}.

Our starting point is based upon ``tropical dynamics'', following \cite{noel2014a,soliman2014a,desoeuvres2022a,kristiansen2026a,portisch2021novel}. In particular, we embed $(\tilde \alpha,\tilde \beta,\tilde \gamma)$ within an $\epsilon$-family:
\begin{align}\eqlab{abc}
 \begin{cases}
  \tilde \alpha =\e^{\alpha\epsilon^{-1}},\\
  \tilde \beta =\e^{\beta\epsilon^{-1}},\\
  \tilde \gamma =\e^{\gamma \epsilon^{-1}},
 \end{cases}\quad (\alpha,\beta,\gamma)\in \mathbb R^3.
\end{align}
In this way, by choosing different $(\alpha,\beta,\gamma)$ we can study different orders of magnitude of $\tilde \alpha$, $\tilde \beta$ and $\tilde \gamma$ by letting $\epsilon\rightarrow 0$. For example, by setting $\beta=1$, $\alpha=\epsilon \alpha_2$, $\gamma=\epsilon \gamma_2$,
we obtain $\tilde \alpha=\mathcal O(1)$, $\tilde \gamma=\mathcal O(1)$, both nonzero, and $\tilde \beta$ large with respect to $\epsilon\rightarrow 0$:
\begin{align*}
 \begin{cases}
  \tilde \alpha =\e^{\alpha_2},\\
  \tilde \beta =\e^{\epsilon^{-1}},\\
  \tilde \gamma =\e^{\gamma_2}.
 \end{cases}
\end{align*}
 In this case, \eqref{mm} is a classical slow-fast system with respect to $\varepsilon: = \e^{-\epsilon^{-1}}\rightarrow 0$ and there is a normally hyperbolic and attracting critical manifold defined by $\tilde u'=0$, i.e.
\begin{align}
 \tilde v = \frac{\tilde u}{\tilde \alpha+\tilde u}.\eqlab{cm0}
\end{align}
Fixing compact sets of $\tilde u$ and $\tilde \alpha>0$, we obtain a slow manifold that is $C^k$-smoothly $\mathcal O(\varepsilon)$-close to the critical manifold for all $0<\varepsilon\ll 1$.
Different choices of $\beta>0$ fixed clearly give equivalent families. Notice, however, if we take $\alpha<0$, then the critical manifold degenerates to $\{\tilde v=1\}$ and $\{\tilde u=0\}$ for $\epsilon\rightarrow 0$. (To see this, notice that $\alpha<0$ corresponds to $\tilde \alpha\rightarrow 0$ as $\epsilon\rightarrow 0$ in \eqref{cm0}.) To study this case with classical methods, one would augment $\epsilon'=0$ and blowup the degenerate set $\tilde v=1$. However, the weights of the blow-up will depend upon the choice of $(\alpha,\beta)$. Different choices of $(\alpha,\beta,\gamma)$ lead to separate degeneracies that through classical methods require separate blowups.

In this paper, we are concerned with all possible orders of magnitudes of $(\tilde \alpha,\tilde \beta,\tilde \gamma)$. However, for simplicity, we will primarily focus on $(\alpha,\beta,\gamma)$ fixed and $\epsilon\rightarrow 0$ in the remainder of the paper. We believe that the remaining cases can be studied using the same methods, but leave this to the interested reader. {We refer to \secref{conclude} for further details.}

We emphasize that $(\alpha,\beta,\gamma)$ is equivalent with any point on the ray $k(\alpha,\beta,\gamma)$, $k>0$, in the limit $\epsilon\to 0$. This is a consequence of \eqref{abc} being invariant with respect to $(\alpha,\beta,\gamma,\epsilon)\mapsto k
(\alpha,\beta,\gamma,\epsilon)$ for any $k>0$. Consequently, one can reduce $(\alpha,\beta,\gamma)$ to the projective space. However, such reduction will not be important to us. 

Now, to study \eqref{mm} with $(\tilde \alpha,\tilde \beta,\tilde \gamma)$ given by \eqref{abc}, we naturally also define $(u,v)$ in the same way
\begin{align}\eqlab{sc0}
 \begin{cases}
  \tilde u =\e^{u\epsilon^{-1}},\\
  \tilde v =\e^{v\epsilon^{-1}}.
 \end{cases}
\end{align}
Notice that by fixing $(u,v)$ in a compact set, this transformation will in principle also allow us to study unbounded domains of $(\tilde u,\tilde v)$-plane as $\epsilon\rightarrow 0$. The system \eqref{mm} has a tropical description in the limit $\epsilon\to 0$ when expressed in the $(u,v)$-variables with $(\tilde \alpha,\tilde \beta,\tilde \gamma)$ given by \eqref{abc}. This is described in details in \cite{kristiansen2026a} \textcolor{black}{for general polynomial systems; see \secref{trop} for a description relating to \eqref{mm}}.

\subsection{Summary of main results}

We will first use the tropical approach of \cite{kristiansen2026a} to obtain a classification of \eqref{mm} with $(\tilde \alpha,\tilde \beta,\tilde \gamma)$ given by \eqref{abc} in the singular limit $\epsilon\to 0$.
This classification is based upon tropical geometry and
leads to three generic regions in the $(\alpha,\beta,\gamma)$-space corresponding to different singular geometries and dynamics in the $(u,v)$-plane. In particular, we identify suitable tropical lines in the $(u,v)$-plane which partially encode the slow dynamics. The full dynamics of \eqref{mm} for $\epsilon>0$ small enough is analyzed by combining blowup procedures with GSPT, see  \cite{jelbart2025switchingmultipletimescalesgeometric,kaklamanos2019a}. Our main result, see \thmref{main}, shows that on a suitably blown up space we find a globally defined and normally hyperbolic critical manifold (which perturbs to a global slow manifold) in all cases except when $\alpha<\gamma<\beta$ (which we refer to as case (iii$2$)). In this case, there is a subset of the blown up space where the unique manifold from infinity only extends by regular flow (as $\epsilon\to 0$).

We anticipate that our approach can be applied to more general multi-scale systems, which can exhibit different types of slow manifolds in different asymptotic regimes of parameters and variables, see e.g. \cite{Gucwa2009783}. In \cite{kristiansen2026a}, the tropical approach was used to classify
asymptotic regimes for the autocatalator.
Finding suitable scalings is challenging in general and systematic approaches are needed. We see our work as a contribution to this important area, see \cite{goeke2014a,noel2014a} for other approaches in this direction.  

\subsection{Outline}
The paper is outlined as follows: In \secref{trop}, we first review the basic concepts from \cite{kristiansen2026a} on tropical dynamics. This leads to a basic classification of the singular limit systems for $\epsilon\to 0$. Then in \secref{main} we present our main result \textcolor{black}{on Problem \ref{problem}}, which we prove in the remaining sections (\secref{casec}--\secref{completing}) of the paper. In \secref{conclude}, we conclude the paper with a discussion section, where we compare our findings with existing work.
%

\section{Tropical dynamics}\seclab{trop}
\textcolor{black}{We consider \eqref{mm} with parameters $(\tilde \alpha,\tilde \beta,\tilde \gamma)$ of the form \eqref{abc} and the change of coordinates $(u,v)\mapsto (\tilde u,\tilde v)$ defined by \eqref{sc0}, repeated here for convenience:}
\begin{align}\eqlab{sc}
 \begin{cases}
  \tilde u =\e^{u\epsilon^{-1}},\\
  \tilde v =\e^{v\epsilon^{-1}}.
 \end{cases}
\end{align}
 Define
$F_i=F_i(u,v;\alpha,\beta,\gamma,\epsilon)$ and $\delta_i$, $i\in \{1,\ldots,6\}$, as follows
\begin{equation}\eqlab{tropauto}
\begin{alignedat}{3}
 F_1(u,v)&:=0,\quad &\delta_1&:=-1,\\
 F_2(u,v)&:=\alpha+v-u,\quad  &\delta_2&:=1,\\
 F_3(u,v)&:=v, \quad &\delta_3&:=1,\\
 F_4(u,v)&:=u-v-\beta, \quad &\delta_4&:=1,\\
 F_5(u,v)&:=\m -\beta,\quad &\delta_5&:=-1,\\
 F_6(u,v)&:=u-\beta,\quad &\delta_6&:=-1,
\end{alignedat}
\end{equation}
with
\begin{align}\eqlab{mepsdef}
 \m:=\epsilon \log \left(\e^{\alpha\epsilon^{-1}}+\e^{\gamma \epsilon^{-1}}\right).
\end{align}
Notice that
\begin{align}\eqlab{meps}
 \m\rightarrow \max(\alpha,\gamma)\quad \mbox{for}\quad \epsilon\rightarrow 0.
\end{align}
Then a simple calculation, based upon
\begin{align}\eqlab{tildexp}
\begin{cases}
  u' = \epsilon \tilde u^{-1}\tilde u',\\
 v' =\epsilon \tilde v^{-1} \tilde v',
 \end{cases}
\end{align}
shows that $(u,v)\mapsto (\tilde u,\tilde v)$ brings \eqref{mm} into the following form:
%
\begin{equation}\eqlab{uv}
\begin{aligned}
 \dot u &=\frac{\sum_{i=1}^3 \delta_i \e^{F_i(u,v)\epsilon^{-1}}}{\sum_{i=1}^6 \e^{F_i(u,v)\epsilon^{-1}}},\\
 \dot v &=\frac{\sum_{i=4}^6 \delta_i \e^{F_i(u,v)\epsilon^{-1}}}{\sum_{i=1}^6 \e^{F_i(u,v)\epsilon^{-1}}},
\end{aligned}
\end{equation}
with dot denoting differentiation with respect to a ``normalized time'', defined implicitly by
\begin{align}\eqlab{time}
 \dot{()} = \frac{\epsilon^{-1} ()'}{{\sum_{i=1}^6 \e^{F_i(u,v)\epsilon^{-1}}}}.
\end{align}

We now consider the tropicalized version of \eqref{uv} \textcolor{black}{obtained from the singular limit $\epsilon \to 0$}.
For this purpose, recall that in the tropical semiring \cite{maclagan2015introduction}, multiplication is defined as addition while addition is defined by $\max$ (it is also possible to use $\min$ but we stick with $\max $ here). As an example, consider the term
\begin{align}\eqlab{term}
\tilde \beta^{-1} (\tilde \alpha+\tilde \gamma)\tilde v,\end{align} and its tropical version $\max(\alpha,\gamma)+v-\beta$; we drop the tildes when going to the tropical sense. The basic mechanism (known as Maslov's dequantization) for the tropicalization is the following:  Substituting \eqref{abc} and \eqref{sc} into \eqref{term}  gives
\begin{align*}
\epsilon \log \left(\tilde \beta^{-1} (\tilde \alpha+\tilde \gamma)\tilde v\right) \to \max(\alpha,\gamma)+v-\beta,
\end{align*}
in line with \eqref{meps},
for $\epsilon\to 0$. This crucial observation has been worked out in detail in \cite{kristiansen2026a} in the context of planar polynomial vector-fields.

Following the concepts and terminology of \cite{kristiansen2026a}, we define
\begin{align*}
 \mathcal U:=\{1,2,3\},\quad\mathcal V:=\{4,5,6\},\quad \mathcal I:=\mathcal U\cup \mathcal V.
\end{align*}
The affine functions $\{F_i\}_{i\in \mathcal U}$ ($\{F_i\}_{i\in \mathcal V}$) in \eqref{tropauto} are therefore (for $\epsilon\rightarrow 0$) the tropicalized versions of the $(\tilde u,\tilde v)$-monomials of $\tilde u^{-1}\tilde u'$ ($\tilde v^{-1}\tilde v'$, resp.); the division by $\tilde u$ and $\tilde v$ is due to \eqref{tildexp}.
%
At the same time, the constants
$\{\delta_i\}_{i=1}^6$ define the signs of the corresponding terms.
We also define the tropical flow vectors:
\begin{align}\eqlab{flowvec}
\md_i = \begin{cases}
                                    (\delta_i,0)^{\mathrm T}, & i\in \mathcal U,\\
                                    (0,\delta_i)^{\mathrm T}, & i\in \mathcal V.
                                   \end{cases}
\end{align}
\textcolor{black}{In this way, we can write \eqref{uv} in the form vector form:}
\begin{align}\eqlab{uv1}
 \begin{pmatrix}
  \dot u\\
  \dot v
 \end{pmatrix} = \frac{\sum_{i\in \mathcal I} \md_i \e^{F_i(u,v)\epsilon^{-1}}}{\sum_{i\in \mathcal I} \e^{F_i(u,v)\epsilon^{-1}}}.
 \end{align}
Finally, we define the tropical curve $\mathcal T$ associated with the monomials $\{F_i\}_{i\in \mathcal I}$ as the set of points $(u,v)$ where $\#\argmax_{i\in \mathcal I} F_i(u,v)\ge 2$, \textcolor{black}{using $\#$ to denote the cardinality of the $\argmax$-set:}
$$\textcolor{black}{\argmax_{i\in \mathcal I} F_i(u,v) := \left\{i \in \mathcal I\,:\,F_i(u,v) = \max_{j\in \mathcal I} F_j(u,v)\right\}}.$$

Now, \eqref{mm} and \eqref{uv}, with $(\tilde \alpha,\tilde \beta,\tilde \gamma)\in (0,\infty)^3$ defined by \eqref{abc} for $(\alpha,\beta,\gamma)\in \mathbb R^3$ and $\epsilon>0$, are by construction smoothly equivalent for any $\epsilon>0$.
However, \eqref{uv} with $(\alpha,\beta,\gamma)$ fixed has a piecewise smooth (PWS) limit as $\epsilon\to 0$, \textcolor{black}{with the discontinuity set being a subset of $\mathcal T$. Indeed, consider any $(u,v)\notin \mathcal T$. Then by definition $\argmax_{i\in \mathcal I} F_i(u,v)\in \mathcal I$ is a singleton and from \eqref{uv1}, we find that}
\begin{align}\eqlab{uv0}
 \begin{pmatrix}
  \dot u\\
  \dot v
 \end{pmatrix} = \frac{   \md_{\argmax_{i\in \mathcal I} F_i(u,v)} +o(1)  }{1+o(1)}\to
 \md_{\argmax_{i\in \mathcal I} F_i(u,v)}\quad \forall\,(u,v)\notin \mathcal T,
\end{align}
as $\epsilon\to 0$, see also \cite[Lemma 3.3]{kristiansen2026a}; this explains why we work with the time defined by \eqref{time}. 

\subsection{Further background}
The reference \cite{kristiansen2026a} defines tropical dynamical systems through a differential inclusion \begin{align*}
\begin{pmatrix}\dot u\\ \dot v\end{pmatrix}\in \Conv(u,v),
\end{align*}
for arbitrary tropical monomials $\{F_i\}_{i\in \mathcal I}$ and arbitrary signs $\{\delta_i\}_{i\in \mathcal I}$, with $\mathcal I=\mathcal U\cup \mathcal V$ where $\mathcal U$ and $\mathcal V$ are disjoint index sets that correspond to monomials from the $u$-direction resp. $v$-direction.
The differential inclusion agrees with the singular limit \eqref{uv0} for $(u,v)\notin \mathcal T$:
\begin{align*}
 \Conv(u,v)=\{\md_{\argmax_{i\in \mathcal I} F_i(u,v)}\}\quad \forall\,(u,v)\notin \mathcal T,
 \end{align*}
On $\mathcal T$, the differential inclusion extends the concepts of crossing and sliding known for Filippov's theory of piecewise smooth systems, see \cite{Bernardo08,filippov1988differential}, with the aim to capture the actual (slow) flow for $\epsilon>0$ small. In \cite{kristiansen2026a} it is shown that there exists a global solution of the differential inclusion through any point whose orbit corresponds to a polygonal curve. This leads the authors to define structural stability of tropical systems through homotopic polygonal curves.  We will not need the full details of \cite{kristiansen2026a}, but we will use the subdivision $\mathcal S$ of the Newton polygon, which also plays an important role in the classification of structural stability in \cite{kristiansen2026a}. To explain the subdivision,
we first note that for a general tropical monomial $F(u,v)=d+nu+mv$, we define $\coeff F: =d\in \mathbb R$ as the associated tropical coefficient of $F$ and $\degree F: = (n,m)\in \mathbb N_{-1}^2$, with $\mathbb N_{-1}:=\{-1,0\}\cup \mathbb N$, as the associated degree.
Consider then the set of points
 $$(\degree F_i,\coeff F_i)\in \mathbb N_{-1}^2\times \mathbb R,\quad i\in \mathcal I,$$  obtained by assigning the tropical coefficients $\coeff F_{i}\in \mathbb R$ of $F_{i}$ as a third dimension (a height). Then we take the upper envelope of the resulting set of points. Finally, the subdivision $\mathcal S$ is obtained by projecting the edges of the resulting convex polyhedron back to $\mathbb R^2$ where $\mathbb N_{-1}^2$ is naturally embedded, see \figref{subdivision}(a). The subdivision is dual to the tropical curve in the following sense (see \cite[Proposition 3.1.6]{maclagan2015introduction}):
Each face of the subdivision corresponds to a {tropical vertex}  $(u,v)\in \mathcal T$ where
\begin{align}\eqlab{tropPcond}
\# \argmax_{i\in \mathcal I} F_{i}(u,v)>2,
\end{align} and the edges of $\mathcal S$ correspond to
the edges of $\mathcal T$ where $$\# \argmax_{i\in \mathcal I} F_{i}(u,v)=2.$$
In particular, the subdivision $\mathcal S$ fixes the tropical curve $\mathcal T$ up to homotopy and the dual edges of $\mathcal S$ and $\mathcal T$ are mutually perpendicular, see \figref{subdivision}(b). The numbers in the boxes at $\degree F_i$, $i\in \{1,\ldots,6\}$, refer to the subscripts $i$ of the monomials. At the same time, the arrows correspond to the associated flow vectors $\md_i$, recall \eqref{flowvec}.
  A subdivision $\mathcal S$ where all of its faces are triangles is called a \textit{triangulation}, and in this case exactly three \rsp{tropical edges} come together at each {tropical vertex}. This is the generic case when all $\coeff F_i$ are free to vary.

  \begin{figure}
  \begin{center}
    \subfigure[]{\includegraphics[width=0.55\textwidth]{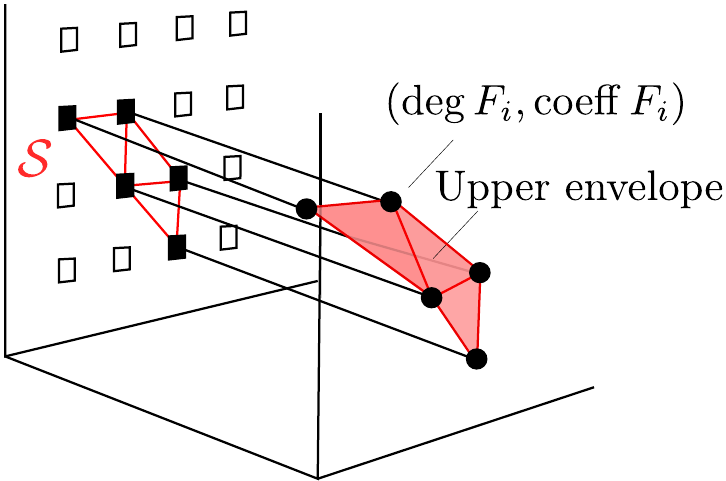}}
        \subfigure[]{\includegraphics[width=0.34\textwidth]{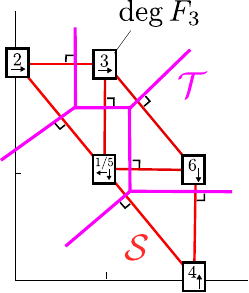}}
    \caption{In (a): The upper envelope of the set of points $(\degree F_{i},\coeff F_{i})$, $i\in \{1,\ldots,6\}$. The subdivision $\mathcal S$ of the Newton polygon is the projection of the envelope. In (b): The subdivision $\mathcal S$ and the dual tropical curve $\mathcal T$. The subdivision fixes $\mathcal T$ up to homotopy with the edges of $\mathcal T$ being perpendicular to the edges of $\mathcal S$. The numbers $1,\ldots,6$ correspond to the subscripts of the monomials in \eqref{tropauto} while the arrows correspond to the associated flow vectors $\md_i$, recall \eqref{flowvec}.}\figlab{subdivision}
        \end{center}
    \end{figure}
\subsection{A first classification}
For \eqref{tropauto}, we have
\begin{equation}\eqlab{tropauto2}
\begin{alignedat}{3}
 \coeff F_1&:=0,\quad &\degree F_1 &:=(0,0),\\
 \coeff F_2&:=\alpha,\quad  &\degree F_2&:=(-1,1),\\
 \coeff F_3&:=0, \quad &\degree F_3&:=(0,1),\\
 \coeff F_4&:=-\beta, \quad &\degree F_4&:=(1,-1),\\
 \coeff F_5&:=\max(\alpha,\gamma) -\beta,\quad &\degree F_5&:=(0,0),\\
 \coeff F_6&:=-\beta,\quad &\degree F_6&:=(1,0).
\end{alignedat}
\end{equation}
We are interested in the possible sub-divisions for all $(\alpha,\beta,\gamma)\in \mathbb R^3$; notice that for \eqref{tropauto} not all $\coeff F_i$ are free to vary.
A basic calculation shows that the points $(\degree F_i,\coeff F_i)$, $i\in \{1,3,4,6\}$, all lie in the same plane given by
\begin{align}
 z =-\beta x,\eqlab{xyz}
 \end{align}
setting $z:=\coeff F,\,(x,y):=\degree F$. From this follows that if $\coeff F_5\le \coeff F_1\Leftrightarrow \max(\alpha,\gamma)\le \beta$, then the subdivision is not a triangulation. Suppose that $\max(\alpha,\gamma)<\beta$. Then $(\degree F_2,\coeff F_2)$ lies ``below'' \eqref{xyz} since for $x=-1$ we have $z=\beta>\coeff F_2=\alpha$ on this plane.  We conclude that we have the subdivision in \figref{subdivision2}(c) for all $\max(\alpha,\gamma)<\beta$.

On the other hand, for $\max(\alpha,\gamma)>\beta$ then the two set of points defined by $(\degree F_i,\coeff F_i)$, $i\in \{3,5,6\}$, and $(\degree F_i,\coeff F_i)$, $i\in \{4,5,6\}$, lie in two different planes. Moreover, the plane containing the former set of points $(\degree F_i,\coeff F_i)$, $i\in \{3,5,6\}$, is given by
\begin{align*}
 z =-\max(\alpha,\gamma)x+(\beta-\max(\alpha,\gamma))(y-1).
\end{align*}
For $(x,y)=\degree F_2=(-1,1)$, we have $z=\max(\alpha,\gamma)\ge \coeff F_2 = \alpha$ on this plane. Consequently, for any $\alpha \ge \gamma$, the points $(\degree F_i,\coeff F_i)$, $i\in \{2,3,5,6\}$, all lie in the same plane. This gives the subdivision in \figref{subdivision2}(b), which also is not a triangulation. Finally, for $\alpha<\gamma$ we have the subdivision in \figref{subdivision2}(a), which is a triangulation.

For $\gamma<\alpha=\beta$, the points $(\degree F_i,\coeff F_i)$, $i\in \{1,2,3,4,6\}$, all lie in the same plane. For $\gamma=\alpha$ then $F_1=F_5$ and the points $(\degree F_i,\coeff F_i)$, $i\in \{1,2,3,4,5,6\}$, all lie in the same plane and the subdivision is trivial. We will not consider these degenerate cases here.

  \begin{figure}
  \begin{center}
  \subfigure[$\beta<\max(\alpha,\gamma)=\gamma$]{\includegraphics[width=0.3\textwidth]{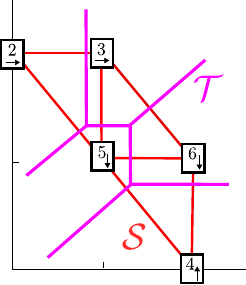}}
        \subfigure[$\beta<\max(\alpha,\gamma)=\alpha$]{\includegraphics[width=0.3\textwidth]{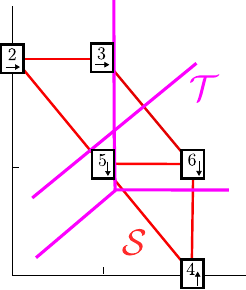}}
    \subfigure[$\max(\alpha,\gamma)<\beta$]{\includegraphics[width=0.3\textwidth]{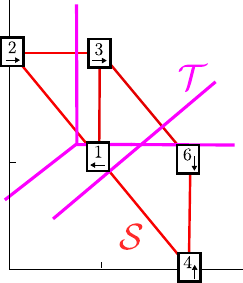}}
    \caption{The three possible subdivisions $\mathcal S$ (in red) associated with the $(\alpha,\beta,\gamma)$-family of the tropical monomials $\{F_i\}_{i\in \mathcal I}$ in \eqref{tropauto}. The associated tropical curves $\mathcal T$ are shown in purple. For further details see the main text.}\figlab{subdivision2}
        \end{center}
    \end{figure}

    In \figref{pp}, we show the tropical phase portraits corresponding to the subdivisions in \figref{subdivision2}. These figures are obtained using the \verb#Matlab# toolbox \verb#Tropical# \verb#Phase# \verb#Plane#  developed by A. Sarantaris in \cite{andreas}. The toolbox can be downloaded from

\begin{center}
\url{https://www.mathworks.com/matlabcentral/fileexchange/132058-tropical-dynamics-toolbox}.
\end{center}

    The tropical curve separates the tropical regions (in different colours with numbers referencing \eqref{tropauto}) where $$\#\argmax_{i\in \mathcal I} F_i(u,v)=1.$$ The dotted lines are crossing lines whereas the full lines are sliding regions. The case $\max(\alpha,\gamma)<\beta$ in (c) actually comes in two versions. The one shown is for $\alpha<\gamma<\beta$ where the tropical system has a downward flow along the edge $\mathcal E_{12}$ corresponding to $F_1(u,v)=F_2(u,v)>F_i$ for all $i\in \{3,4,5,6\}$. This flow along $\mathcal E_{12}$ is called tropical sliding flow (proxy for slow flow) in \cite{kristiansen2026a} and is due to the fact that $F_5$  (with $\delta_5=-1$) dominates the $v$-monomials along $\mathcal E_{12}$.  The red dotted lines indicate the tropical curves obtained from the monomials $F_i$, $i\in \mathcal V$, in the $v$-direction, whereas the blue dotted lines indicate the tropical curves obtained from the monomials $F_i$, $i\in \mathcal U$, in the $u$-direction; these dotted lines are only shown in regions where they differ from the tropical curve. For $\gamma\le \alpha<\beta$, the red and blue lines coincide in (c) and consequently there is no tropical sliding flow along any of the sliding segments (both directions are equilibrated along $\mathcal E_{12}$). Notice that the segment $\mathcal E_{13}$ corresponding to $F_1(u,v)=F_3(u,v)>F_i(u,v)$ for all $i\in \{2,4,5,6\}$ is of so-called tangential nullcline sliding, see \cite[Section 3.2]{kristiansen2026a}. In particular, $\md_1=(-1,0)^{\mathrm T}$ and $\md_3=(1,0)^{\mathrm T}$ are both tangent to $\mathcal E_{13}$; as a tropical dynamical system one can therefore both move to the left and to the right.

    In \figref{pp}(b) there is also no tropical sliding flow along the sliding segments $\mathcal E_{45}$ and $\mathcal E_{46}$, corresponding to $F_4(u,v)=F_5(u,v)>F_i(u,v)$ for all $i\in\{1,2,3,6\}$ resp. $F_4(u,v)=F_6(u,v)>F_i(u,v)$ for all $i\in\{1,2,3,5\}$. In contrast, in the case \figref{pp}(a) there is a downward flow along $\mathcal E_{45}$ since here $F_1(u,v)$ with $\delta_1=-1$ dominates the $u$-monomials: $F_1(u,v)>F_i(u,v)$ for all $i\in \mathcal U\setminus \{1\}=\{2,3\}$ and all $(u,v)\in \mathcal E_{45}$ (notice the blue dotted lines indicating the tropical curve associated with $\{F_{i}\}_{i\in \mathcal U}$).

    In conclusion, we obtain four different equivalence classes of structurally stable systems: (i): $\beta<\max(\alpha,\gamma)=\gamma$, (ii): $\beta<\max(\alpha,\gamma)=\alpha$ and finally (iii$1$): $\gamma<\alpha<\beta$, (iii$2$): $\alpha <\gamma<\beta$. We refer to \cite[Section 7]{kristiansen2026a} for details on tropical equivalence and structural stability. However, it is important to notice that from this simple classification of the tropical systems, we immediately have that slow manifolds for any $0<\epsilon\ll 1$ can only exist near $\mathcal E_{45}$ or $\mathcal E_{46}$ in \figref{pp}(a) and (b) (corresponding to cases (i) and (ii)) and near $\mathcal E_{12}$, $\mathcal E_{13}$ or $\mathcal E_{46}$ in \figref{pp}(c) (corresponding to case (iii)). Since the segment $\mathcal E_{13}$ is of tangential type, its attractivity is not visible on the tropical level. It is only established after blowup. The remaining segments $\mathcal E_{12}$, $\mathcal E_{45}$, $\mathcal E_{46}$ are called tropical nullclines in \cite{kristiansen2026a}.

     \begin{figure}
  \begin{center}
                  \subfigure[Case (i): $\beta<\max(\alpha,\gamma)=\gamma$]{\includegraphics[width=0.475\textwidth]{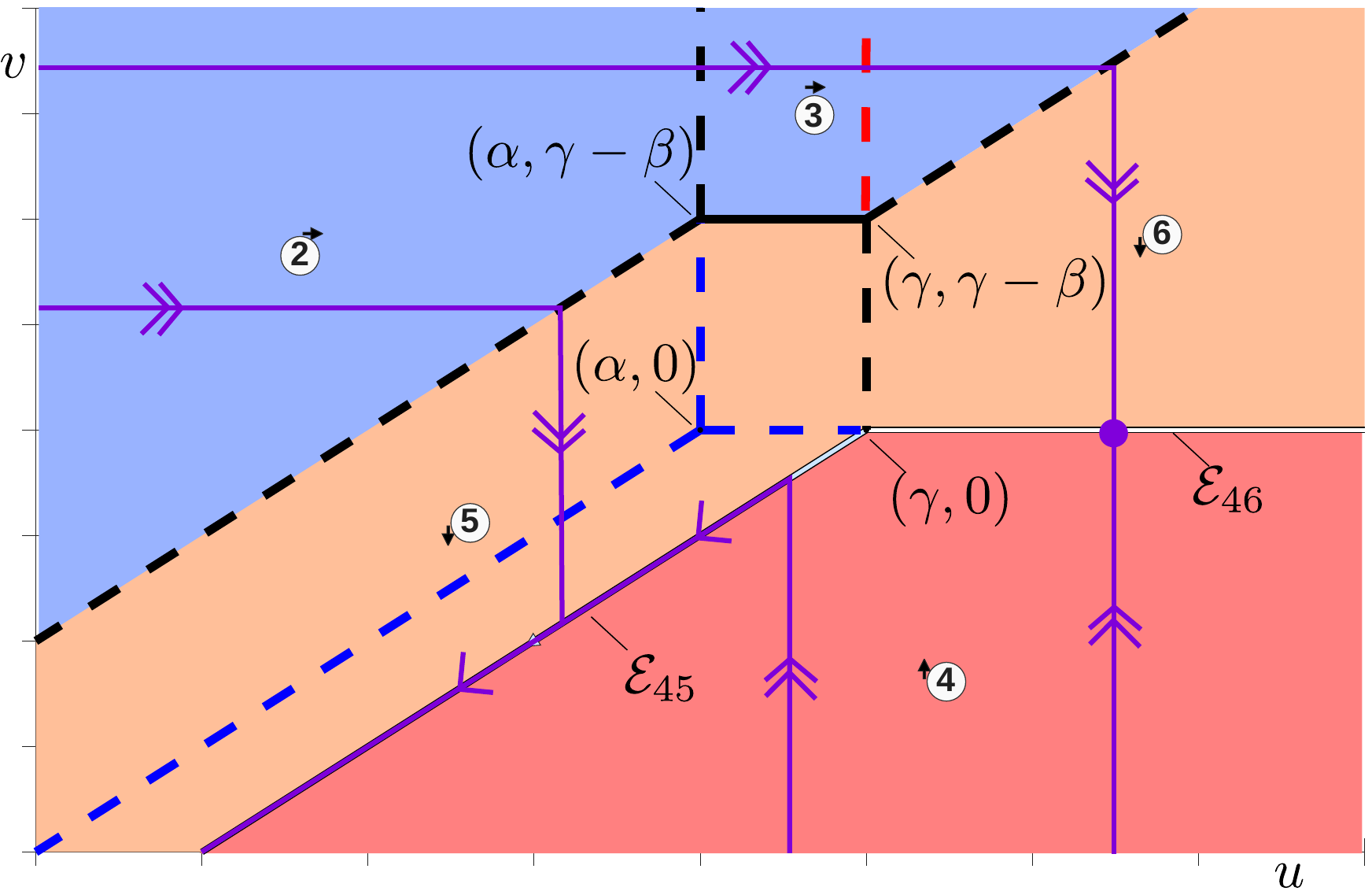}}

        \subfigure[Case (ii): $\beta<\max(\alpha,\gamma)=\alpha$]{\includegraphics[width=0.475\textwidth]{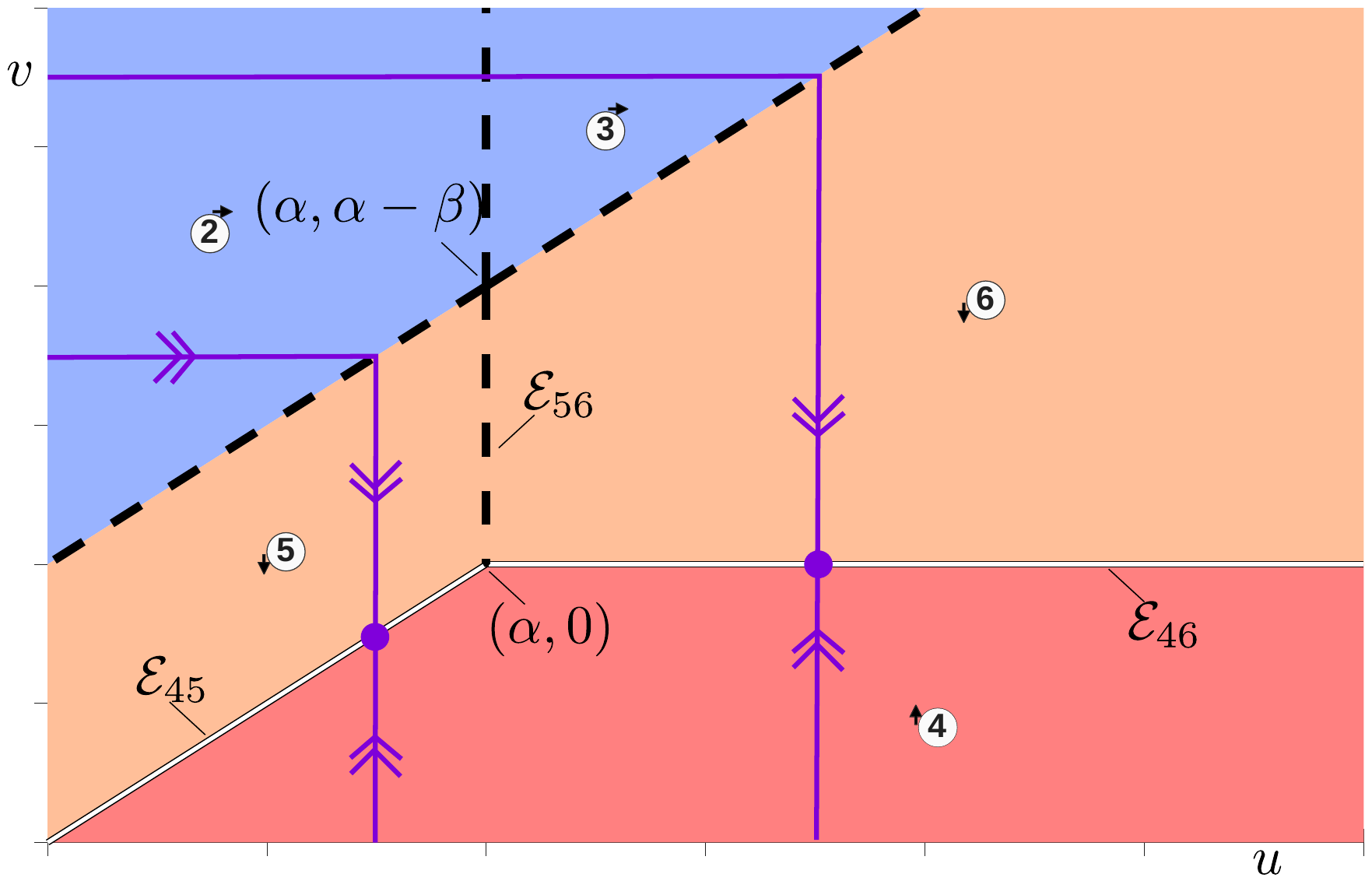}}
    \subfigure[Case (iii): $\max(\alpha,\gamma)<\beta$]{\includegraphics[width=0.475\textwidth]{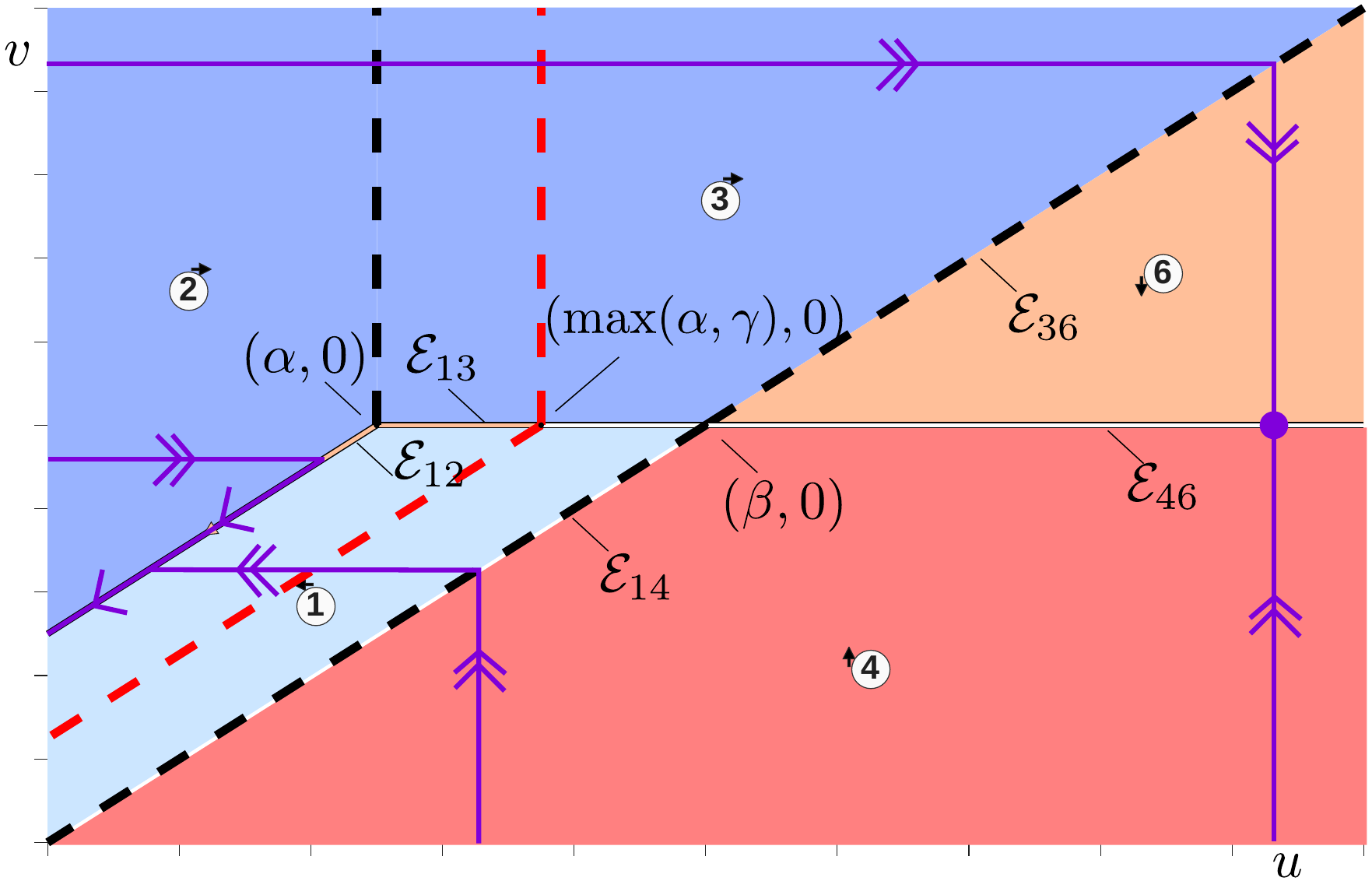}}
    \caption{The tropical phase portraits associated with the subdivisions in \figref{subdivision2}. In (c), we only show the case $\alpha<\gamma<\beta$. The case $\gamma\le \alpha<\beta$ is similar but here there is no tropical sliding flow along $\mathcal E_{12}$.}\figlab{pp}
        \end{center}
    \end{figure}

    \section{Main result}\seclab{main}
    Following the classification above, we restrict attention to the following ``generic'' cases:
    \begin{enumerate}
          \item[(i)] $\max(\alpha,\beta)<\gamma$.
     \item[(ii)] $\max(\beta,\gamma)<\alpha$.
     \item[(iii)] $\max(\alpha,\gamma)<\beta$ and $\alpha\ne \gamma$. In particular, we refer to $\gamma<\alpha<\beta$ as subcase (iii$1$) and $\alpha<\gamma<\beta$ as subcase (iii$2$).
    \end{enumerate}
The remaining cases are ``special cases'' where at least one of the inequalities in (i), (ii) or (iii) is replaced by an equality.  These cases include the classical QSSA regimes (where existence of (local) slow manifolds is known by classical methods). We discuss them in \secref{conclude}. 

Throughout, we use the following definition.
\begin{definition}\defnlab{global}
 Let $(\alpha,\beta,\gamma)$ be fixed, recall \eqref{abc}. We say that a family of one-dimensional invariant manifolds $\widetilde S_\epsilon$ of \eqref{mm}, $0<\epsilon\ll 1$, is a \textnormal{global slow manifold} if $\widetilde S_\epsilon$ connects $\tilde u=\infty$ with the equilibrium $(\tilde u,\tilde v)=(0,0)$ for all  $0<\epsilon \ll 1$, and if there is a sequence of blowups and desingularizations of the extended system defined by \eqref{uv} and $\epsilon'=0$ such that, in each chart of the blown space, $\widetilde S_\epsilon$ converges as $\epsilon\to 0$, uniformly on compact sets, to a normally hyperbolic critical manifold $S_0$, and $S_0$ (taken over all charts) is a connected curve on the blown up space.
\end{definition}

We now summarize our main result, \textcolor{black}{which addresses Problem \ref{problem}}, as follows.
    \begin{theorem}\thmlab{main}
    Consider $(\alpha,\beta,\gamma)$ fixed.
    Suppose first \textnormal{(i)}, \textnormal{(ii)} or \textnormal{(iii$1$)}. Then there is a global slow manifold in the sense of \defnref{global}. 
    On the other hand, in the case \textnormal{(iii$2$)} there is a sequence of blowups and desingularizations such that there are two disjoint normally hyperbolic critical manifolds (graphs over separate/disjoint domains) that are connected  by regular flow. Hence a global slow manifold does not exist in the case \textnormal{(iii$2$)}. 
    \end{theorem}

We will present more details of \thmref{main} below. In particular, we will compute asymptotic expansions of the slow manifolds in certain local charts. Moreover, in the case (iii$2$), we will show that the unique invariant manifold from infinity is a slow manifold all the way down to $u\gtrsim\frac{\beta+\gamma}{2}$, $v=\mathcal O(\epsilon)$ and all $0<\epsilon\ll 1$, which upon extension of the flow enters a region of phase space where the flow is regular (no critical manifold). This region is defined by
\begin{align}\eqlab{u2v21def}
\begin{cases}
 u= \frac{\beta+\gamma}{2}+\epsilon \log x,\\
 v = \epsilon \e^{\frac{\gamma-\beta}{2}\epsilon^{-1}} y,
 \end{cases}
\end{align}
with $x,y=\mathcal O(1)$, $x>0$, $y\in \mathbb R$, as $\epsilon\to 0$, or in terms of the original parameters and variables (using \eqref{abc} and \eqref{sc}) as:
\begin{align}\eqlab{tildeuvxy}
\begin{cases}
 \tilde u= \sqrt{\tilde \beta \tilde \gamma} x,\\
 \tilde v = \e^{\sqrt{\frac{\tilde \gamma}{\tilde \beta}} y},
 \end{cases}
\end{align}
with $\alpha<\gamma<\beta$ and $\epsilon\to 0$ translating into $$\tilde \alpha\ll \tilde \gamma\ll \tilde \beta.$$ {Consequently,  we also have from the $\tilde v$-equation that}
\begin{align*}
 \tilde v = 1+\sqrt{\frac{\tilde \gamma}{\tilde \beta}}y(1+o(1)),
\end{align*}
as $\epsilon\to 0$. By viewing \eqref{tildeuvxy} as a change of coordinates $(x,y)\mapsto (\tilde u,\tilde v)$, we
obtain the following system from \eqref{mm} within compact subsets of $x>0$, $y\in \mathbb R$, and $\epsilon\to 0$, upon appropriate desingularization:
\begin{equation}\eqlab{xyeqn}
\begin{aligned}
 x' &=xy,\\
 y' &= -xy-1.
\end{aligned}
\end{equation}
{This system has no singularities, only regular flow.}
The phase portrait of \eqref{xyeqn} is available in \figref{phaseportrait} below. Here the unique manifold from infinity enters as $y=-x^{-1}(1+\mathcal O(x^{-1}))$ for $\epsilon\to 0$ (in red in \figref{phaseportrait}), see further details in \secref{casea1}.

   The results of the blowups and desingularizations for cases \textnormal{(i)} and \textnormal{(ii)} are illustrated in \figref{ppbu}(a) and \figref{ppbu}(b). The normally hyperbolic and attracting critical manifold is in red. Notice that we here use blowup and desingularization (following \cite{kristiansen2018a,krihog2,kristiansen2020a,uldall2021a}) as a method to resolve the lack of smoothness for $\epsilon\to 0$. A similar procedure for case (iii) leads to the diagram in \figref{ppau}, but here we find a critical set on the blowup of $\mathcal E_{13}$ that remains unresolved (as a nonhyperbolic critical set). The degeneracy is here related to flat terms that cannot be resolved by the classical blowup procedure (now interpreted as a method to gain hyperbolicity, as in \cite{dumortier_1996,krupa_extending_2001}).

    We provide the details of the blowup analysis that relate to slow manifolds in the cases (i) and (iii) in \secref{casec} resp. \secref{casea} below; the case (ii) is similar to (i) and the details are therefore left out (see also \secref{completing} below). To be more precise, we present the details of the blowup in a neighborhood of the point $(u,v,\epsilon)=(\gamma,0,0)$ in the case (i) and in a neighborhood of the point $(u,v,\epsilon)=(\beta,0,0)$ in the case (iii). The details of $(u,v,\epsilon)=(\alpha,0,0)$ in the case (iii) is similar to $(u,v,\epsilon)=(\gamma,0,0)$ in the case (i) and therefore also left out (see also \secref{completing} below). We similarly leave out details of the separate blowups that lead to the diagrams in \figref{ppbu} and \figref{ppau} since these are not relevant for slow manifolds and for the proof of \thmref{main}.

    To resolve the details of $\mathcal E_{13}$ in the case (iii), we use the approach of \cite{kristiansen2017a} to deal with degeneracies stemming from flat terms, see \secref{flat}--\secref{casea1}. Finally, we complete the proof in \secref{completing}.

     \begin{figure}
  \begin{center}
          \subfigure[$\max(\alpha,\beta)<\gamma$]{\includegraphics[width=0.75\textwidth]{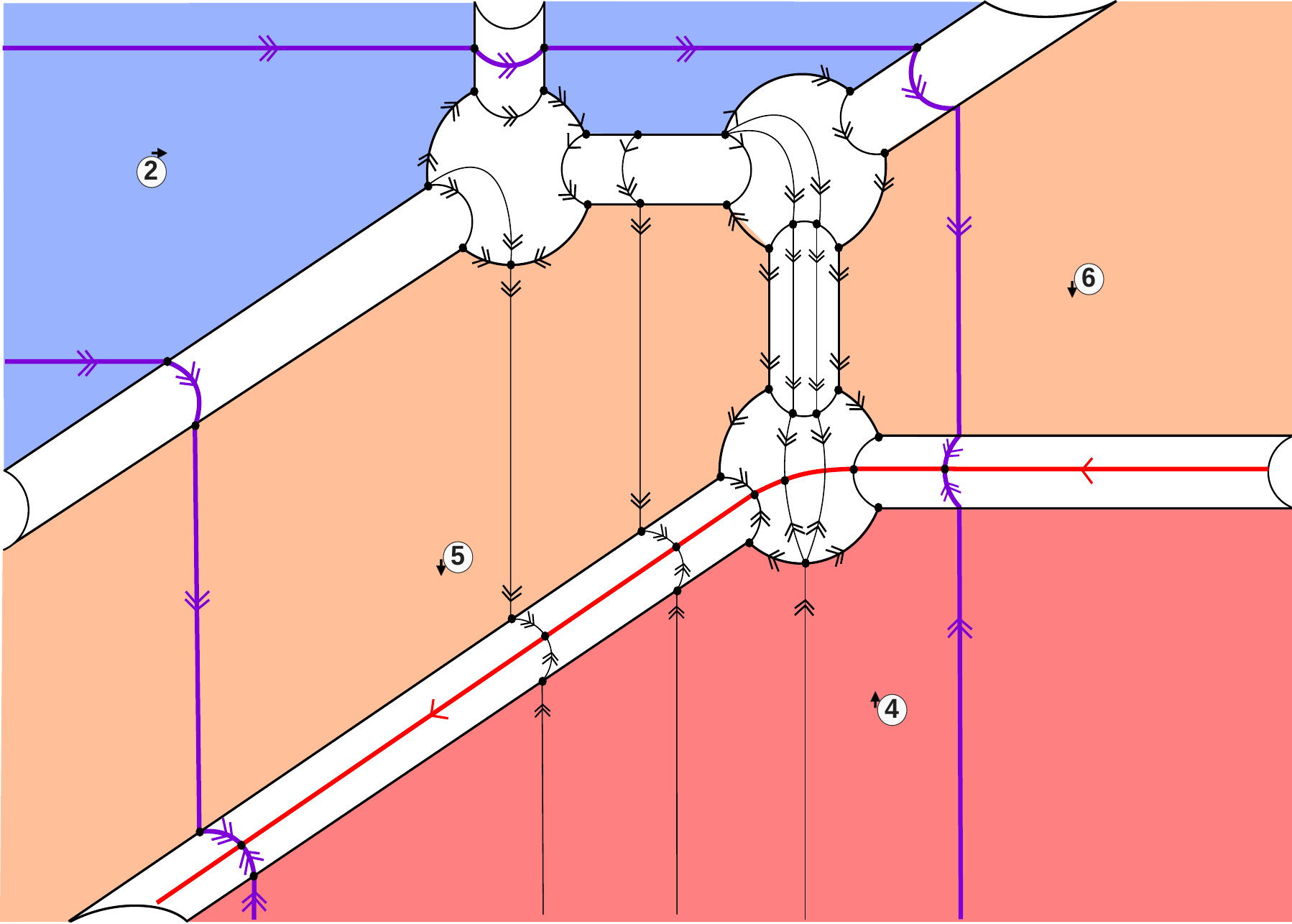}}

    \subfigure[$\max(\beta,\gamma)<\alpha$]{\includegraphics[width=0.75\textwidth]{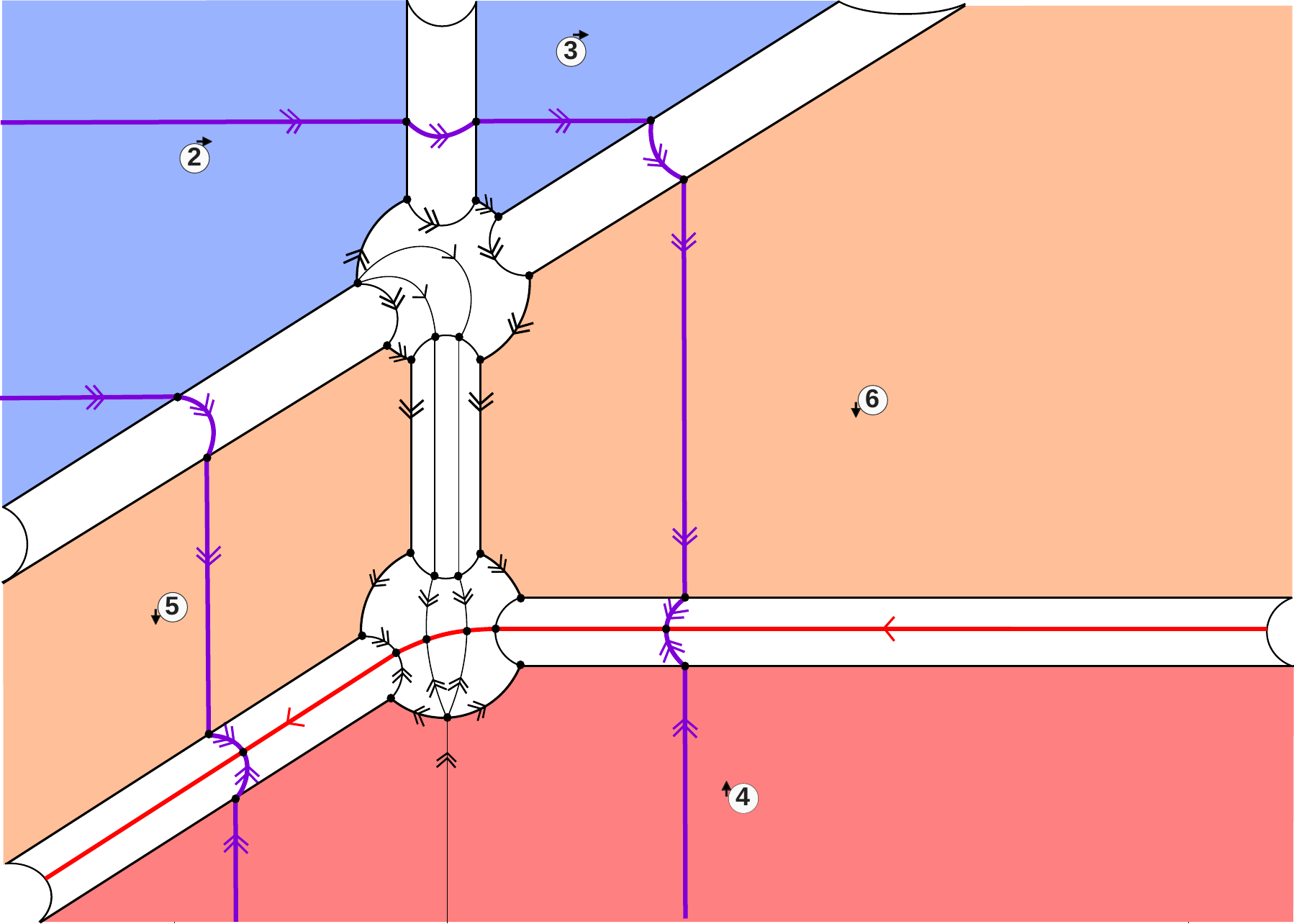}}
%
    \caption{Dynamics on the blowup space for cases (i) and (ii). In both cases, we find a globally defined normally hyperbolic critical manifold on the blown up space. Notice that the third dimension coming out of the page corresponds to $\epsilon$. }\figlab{ppbu}
        \end{center}
    \end{figure}

      \begin{figure}
  \begin{center}
    {\includegraphics[width=0.75\textwidth]{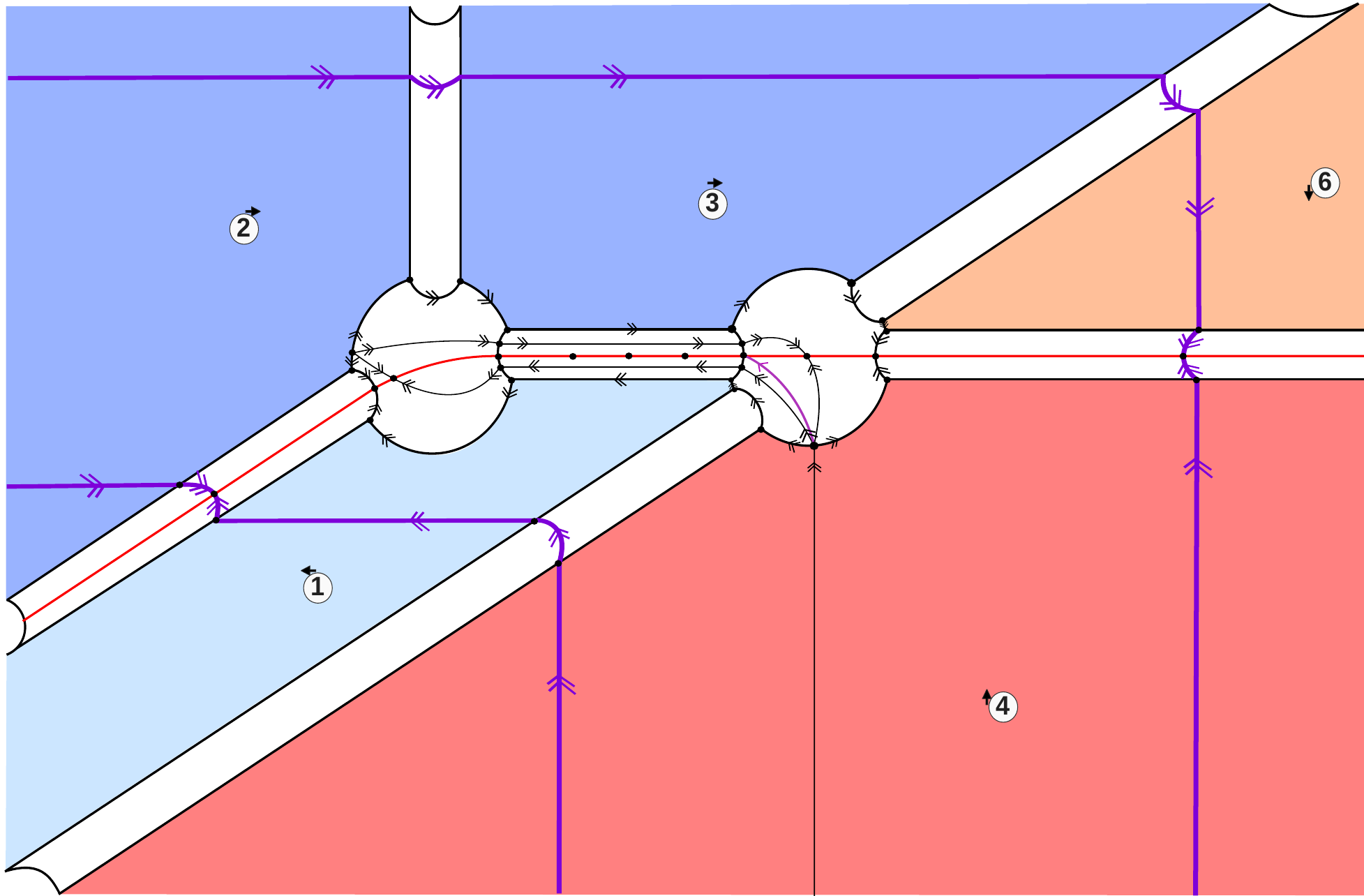}}
%
    \caption{Dynamics on the blowup space for case (iii): $\max(\alpha,\gamma)<\beta$, $\alpha\ne \gamma$. Here we find a critical manifold on the blowup of $\mathcal E_{13}$ (the cylinder between the two spheres) that is degenerate. We resolve this by a separate blowup in  \secref{flat}--\secref{casea1} below, see also \secref{completing}. Notice that the third dimension coming out of the page corresponds to $\epsilon$.  }\figlab{ppau}
        \end{center}
    \end{figure}


\section{Case (i)} \seclab{casec}
In case (i), we have $\max(\alpha,\beta)<\gamma$. The corresponding tropical dynamics is illustrated in \figref{pp}(a). We consider  \eqref{uv} with $\dot \epsilon=0$ augmented in the following form
\begin{equation}\eqlab{uvext}
\begin{aligned}
 u' &=\epsilon \frac{\sum_{i=1}^3 \delta_i \e^{F_i(u,v)\epsilon^{-1}}}{\sum_{i=1}^6 \e^{F_i(u,v)\epsilon^{-1}}},\\
 v' &=\epsilon \frac{\sum_{i=4}^6 \delta_i \e^{F_i(u,v)\epsilon^{-1}}}{\sum_{i=1}^6 \e^{F_i(u,v)\epsilon^{-1}}},\\
 \epsilon' &= 0,
\end{aligned}
\end{equation}
after multiplication of the right hand side by $\epsilon$,
and then blowup $(u,v,\epsilon)=(\gamma,0,0)$ using the following spherical blowup transformation:
\begin{align}\eqlab{bucasec1}
 r\ge 0,\,(\bar u,\bar v,\bar \epsilon)\in S^2\mapsto \begin{cases}
                                                       u =\m+r\bar u,\\
                                                       v =r\bar v,\\
                                                       \epsilon =r\bar \epsilon.
                                                      \end{cases}
\end{align}
Notice that by \eqref{mepsdef}
\begin{align*}
\m = \gamma+\epsilon \log (1+\e^{(\alpha-\gamma)\epsilon^{-1}}) = \gamma+\mathcal O(\epsilon \e^{(\alpha-\gamma)\epsilon^{-1}} ),
\end{align*}
as $\epsilon\rightarrow 0$ by assumption. We work with the following charts
\begin{align}
 \bar u&=1:\quad \begin{cases}
                 u =\m+r_1,\\
                 v =r_1 v_1,\\
                 \epsilon =r_1\epsilon_1,
                \end{cases}\eqlab{c1}\\
                \bar \epsilon&=1:\quad \begin{cases}
                 u =\m+r_2u_2,\\
                 v =r_2 v_2,\\
                 \epsilon =r_2,
                \end{cases}\eqlab{c2}\\
                \bar u&=-1:\quad \begin{cases}
                 u =\m-r_3,\\
                 v =r_3 v_3,\\
                 \epsilon =r_3\epsilon_3.
                \end{cases}\eqlab{c3}
\end{align}
As usual in blowup for smooth systems approaching nonsmooth ones, see \cite{jelbart2021a,uldall2021a,kristiansen2020a}, we achieve the desingularization upon dividing the pull-back vector-field by $\bar \epsilon$; we use $\epsilon_1$, $1$ resp. $\epsilon_3$ in the three charts.

\begin{remark}
 Notice that the $\bar \epsilon=1$-chart corresponds to
 \begin{align*}
 \begin{cases}
  \tilde u = (\tilde \alpha+\tilde \gamma) \e^{u_2},\\
\tilde v =  \e^{v_2}.
\end{cases}
 \end{align*}
 with respect to the original parameters and variables (recall \eqref{abc} and \eqref{sc}) in \eqref{mm}. Here we have used the definition of $\m$ in \eqref{mepsdef} to write
 \begin{align}\eqlab{mepsdef2}
  \e^{\m \epsilon^{-1}} = \e^{\alpha \epsilon^{-1}}+\e^{\gamma \epsilon^{-1}}.
 \end{align}

\end{remark}

The changes of coordinates are defined by the expressions
\begin{align*}
 \begin{cases}
  r_1 =r_2 u_2,\\
  v_1 =u_2^{-1} v_2,\\
  \epsilon_1 =u_2^{-1},
 \end{cases}\quad
\begin{cases}
  r_3 =-r_2 u_2,\\
  v_3 =(-u_2)^{-1} v_2,\\
  \epsilon_3 =(-u_2)^{-1},
 \end{cases}
\end{align*}
for $u_2>0$ resp. $u_2<0$.

In the charts $\bar u=\pm 1$, we find additional degeneracies due to the nonsmoothness along $\mathcal E_{45}$ and $\mathcal E_{46}$, see \figref{pp}(a), that we (as in \cite{jelbart2025switchingmultipletimescalesgeometric,kaklamanos2019a}) resolve using additional cylindrical blowups. For example, in $\bar u=1$, we blow up the degenerate line $r_1\ge 0$, $(v_1,\epsilon_1)=(0,0)$, as follows
 \begin{align}\eqlab{bu1}
 \rho \ge 0,\,(\bar v_1,\bar \epsilon_1)\in S^1\mapsto \begin{cases}
                                                       v_1 =\rho \bar v_1,\\
                                                       \epsilon_1 =\rho \bar \epsilon_1,
                                                      \end{cases}
\end{align}
leaving $r_1$ untouched.
Similarly, in $\bar u=-1$, we blow up the line $r_3\ge 0$, $(v_3,\epsilon_3)=(-1,0)$ as follows:
\begin{align}\eqlab{bu3}
 \rho \ge 0,\,(\bar v_3,\bar \epsilon_3)\in S^1\mapsto \begin{cases}
                                                       v_3 =-1+\rho \bar v_3,\\
                                                       \epsilon_3 =\rho \bar \epsilon_3,
                                                      \end{cases}
\end{align}
leaving $r_3$ untouched.
For simplicity, we will restrict attention to the following charts
\begin{equation}\eqlab{ci2}
\begin{aligned}
                 \bar \epsilon_i&=1:\quad \begin{cases}
                 v_i =-\delta_{i3}+\rho_{i2}v_{i2},\\
                 \epsilon_i =\rho_{i2},
                \end{cases}
\end{aligned}
\end{equation}
for $i\in \{1,3\}$, resp. Here $\delta_{ij}$ denotes Kronecker's delta. Notice that upon composing this with \eqref{c1} resp. \eqref{c3}, we find that
\begin{align*}
 v = \delta_{i3}(u-\m)+\epsilon v_{i2},
\end{align*}
for $i\in \{1,3\}$.
Hence $v_2=v_{12}$, see \eqref{c1}. 

%
%
%

We illustrate the result of the blowup sequence in \figref{ppcbulocal}. Notice in comparison with \figref{ppbu}(b), that our analysis here only limits to the sets relevant to slow manifolds.  The details of resolving the remaining degeneracies are similar and therefore left out for simplicity.


 \begin{figure}
  \begin{center}
    {\includegraphics[width=0.875\textwidth]{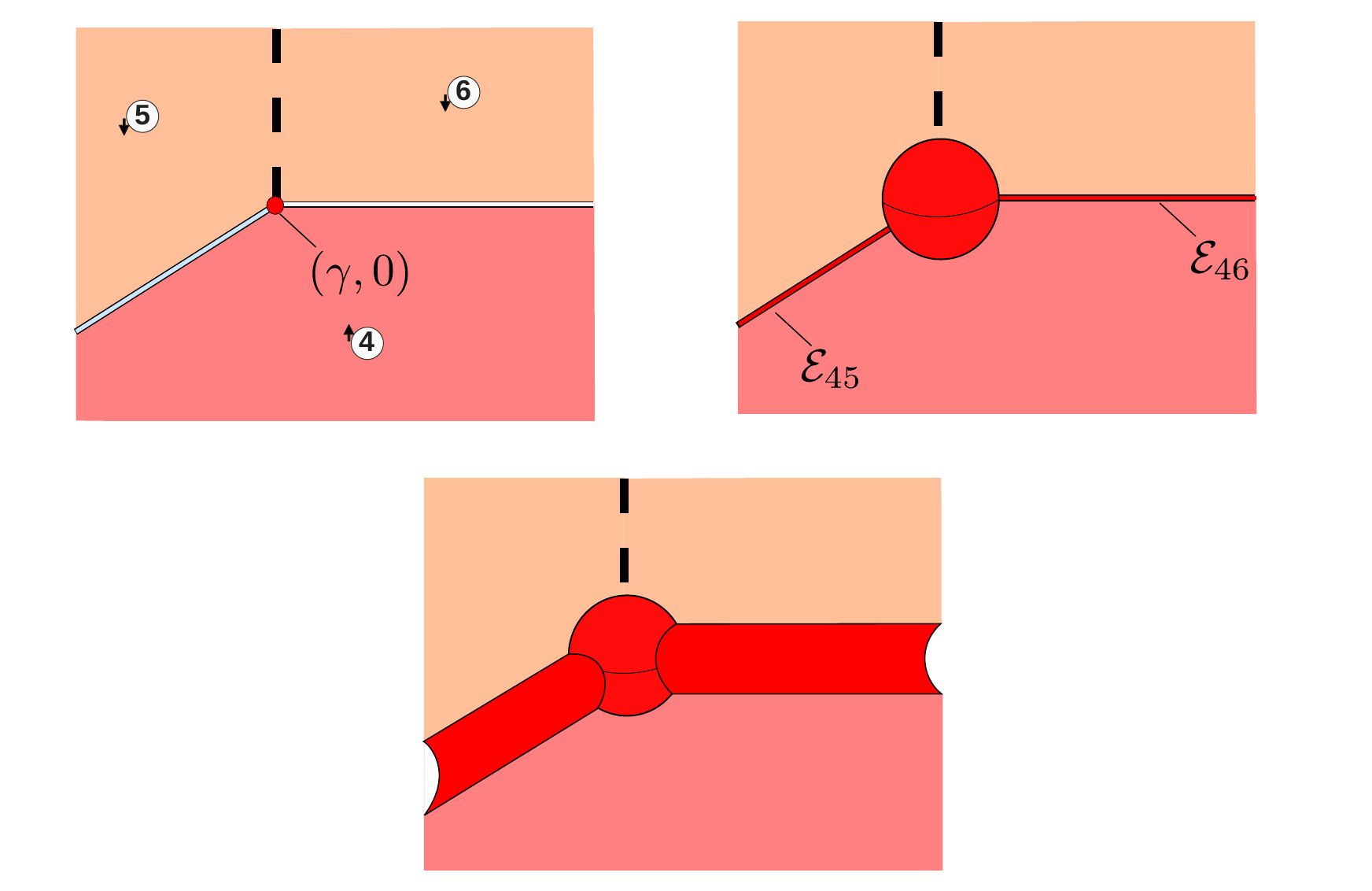}}
%
    \caption{Illustration of the blowup sequence defined by \eqref{bucasec1} and \eqref{bu1}, \eqref{bu3}, going clockwise starting from the top left. First we blow up $(u,v,\epsilon)=(\gamma,0,0)$, indicated in the top left by a red dot in the $(u,v)$-plane, to a sphere (top right). Notice that the third dimension coming out of the page corresponds to $\epsilon$.  Subsequently, we blow up the lines $\mathcal E_{45} \times \{0\}$ and $\mathcal E_{46}\times \{0\}$ to cylinders (in red). The details of the blowup of the black dotted line (corresponding to $\mathcal E_{56}\times \{0\}$, see \figref{pp}(a)) to a third cylinder are left out for simplicity as this does not relate to slow manifolds. }\figlab{ppcbulocal}
        \end{center}
    \end{figure}

\subsection{Analysis in the $\bar \epsilon=1$-chart}
We first consider the $\bar \epsilon=1$-chart. We fix a sufficiently large compact domain $D_2\subset \mathbb R^2$ and consider the pull-back of the vector-field \eqref{uvext} by the map $(u_2,v_2,r_2)\mapsto (u,v,\epsilon)$ defined by \eqref{c2}. This leads to the following equations
\begin{equation}\eqlab{uv2}
\begin{aligned}
 u_2' &=\e^{-(\m-\beta)\epsilon^{-1}}\left[-1+\e^{v_2}+\e^{(\alpha-\m)\epsilon^{-1}}\e^{v_2-u_2}\right],\\
 v_2' &= \e^{u_2-v_2}-1-\e^{u_2},
\end{aligned}
\end{equation}
and $r_2'=\epsilon'=0$, after multiplication of the right hand side by $$\e^{-(\m-\beta)\epsilon^{-1}}\sum_{i=1}^6 \e^{F_i(\m+\epsilon u_2,\epsilon v_2)\epsilon^{-1}},$$ which is regular and nonzero for $(u_2,v_2)\in D_2$ as $\epsilon\rightarrow 0$.
For $\epsilon\to 0$, we obtain
\begin{equation}\nonumber
\begin{aligned}
 u_2' &=0,\\
 v_2' &= \e^{u_2-v_2}-1-\e^{u_2},
\end{aligned}
\end{equation}
since $\gamma>\beta$.
\begin{lemma}
 The system \eqref{uv2} is slow-fast with respect to $\epsilon\to 0$ with a normally hyperbolic and attracting critical manifold $C_2$ defined by
 \begin{align*}
  \e^{v_2} = \frac{\e^{u_2}}{\e^{u_2}+1}.
 \end{align*}
We have the following expansion of an associated slow manifold $S_{2,\epsilon}$ with $S_{2,0}=C_2\cap D_2$:
\begin{align*}
  \e^{v_2} = \frac{\e^{u_2}}{\e^{u_2}+1}+\e^{-(\m-\beta)\epsilon^{-1}}\left(\frac{e^{u_2}}{(1+e^{u_2})^4}+\mathcal O\left(\e^{-(\m-\beta)\epsilon^{-1}},\e^{(\alpha-\m)\epsilon^{-1}}\right)\right) 
\end{align*}
for all $0\le \epsilon<\epsilon_0$, with $\epsilon_0>0$ small enough, and where the remainder is smooth with respect $(\e^{u_2},\e^{-(\m-\beta)\epsilon^{-1}},\e^{(\alpha-\m)\epsilon^{-1}})$.
\end{lemma}
\begin{proof}
 The first statement is obvious. Now, for the expansion of the slow-manifold, we define $\tilde u_2=\e^{u_2}$, $\tilde v_2=\e^{v_2}$ so that
   \begin{equation}\eqlab{x2y2}
\begin{aligned}
 \tilde u_2' &=\e^{-(\m-\beta)\epsilon^{-1}}\left[-\tilde u_2+\tilde v_2\tilde u_2+\e^{(\alpha-\m)\epsilon^{-1}}\tilde v_2\right],\\
 \tilde v_2' &= \tilde u_2-\tilde v_2(1+\tilde u_2).
\end{aligned}
\end{equation}
The result then follows from a simple calculation.
\end{proof}
On $S_{2,\epsilon}$ we have the following slow flow
\begin{align*}
u_2'&=\e^{-(\gamma-\beta)\epsilon^{-1}}\left[-\frac{\e^{u_2}}{1+\e^{u_2}}+\mathcal O\left(\e^{-(\gamma-\beta)\epsilon^{-1}},\e^{(\alpha-\gamma)\epsilon^{-1}}\right)\right]
\end{align*}
with the remainder being smooth with respect $(\e^{u_2},\e^{-(\gamma-\beta)\epsilon^{-1}},\e^{(\alpha-\gamma)\epsilon^{-1}})$. Here we have used \eqref{mepsdef2}. Hence $u_2$ is decreasing for all $0<\epsilon\ll 1$.
\subsection{Analysis in the $\bar u=1$-chart}
Next, in the $\bar u=1$-chart, we fix a compact set $D_1\subset \mathbb R^3$ defined by $(r_1,v_1,\epsilon_1)\in [0,r_{10})\times I\times [0,\epsilon_{10}]$, with $r_{10},\epsilon_{10}>0$ both sufficiently small and $I\subset \mathbb R$ a compact interval containing $0$. The pull-back of the vector-field \eqref{uvext} by the map $(r_1,v_1,\epsilon_1)\mapsto (u,v,\epsilon)$ defined by \eqref{c1}, then takes the following form
\begin{equation}\eqlab{c1eqn}
\begin{aligned}
r_1 '&= \frac{\e^{-(\m-\beta)\epsilon^{-1}}}{1+\e^{\epsilon_1^{-1}}(1+\e^{-\epsilon_1^{-1} v_1})}\left[-1+\e^{(\alpha-\m)\epsilon^{-1}} \e^{-\epsilon_1^{-1}}\e^{\epsilon_1^{-1} v_1}+\e^{\epsilon_1^{-1} v_1} \right]r_1,\\
 v_1'&=\frac{(\e^{-\epsilon_1^{-1} v_1}-1)-\e^{-\epsilon_1^{-1}} }{1+\e^{-\epsilon_1^{-1} v_1}+\e^{-\epsilon_1^{-1}}}\\
 &-\frac{\e^{-(\m-\beta)\epsilon^{-1}}}{1+\e^{\epsilon_1^{-1}}(1+\e^{-\epsilon_1^{-1} v_1})}\left[-1+\e^{(\alpha-\m)\epsilon^{-1}} \e^{-\epsilon_1^{-1}}\e^{\epsilon_1^{-1} v_1}+\e^{\epsilon_1^{-1} v_1} \right]v_1,\\
\epsilon_1 '&= -\frac{\e^{-(\m-\beta)\epsilon^{-1}}}{1+\e^{\epsilon_1^{-1}}(1+\e^{-\epsilon_1^{-1} v_1})}\left[-1+\e^{(\alpha-\m)\epsilon^{-1}} \e^{-\epsilon_1^{-1}}\e^{\epsilon_1^{-1} v_1}+\e^{\epsilon_1^{-1} v_1} \right]\epsilon_1.
\end{aligned}
 \end{equation}
Here we have as discussed above used a desingularization corresponding to division of the right hand side by $\epsilon_1$. Moreover, we have for convenience also multiplied the right hand side by
\begin{align*}
 \frac{\e^{-(\m-\beta)\epsilon^{-1}}\sum_{i=1}^6 e^{F_i(\m+r_1,r_1 v_1)\epsilon^{-1}} }{1+\e^{\epsilon_1^{-1}}(1+\e^{-\epsilon_1^{-1} v_1})},
\end{align*}
which is regular and nonzero on $D_1$. We have left $\epsilon=r_1\epsilon_1$ as a constant where it makes sense. The system clearly extends smoothly to $\{r_1=0\}$ and $\{\epsilon_1=0\}$ provided that $v_1\ne 0$. In fact, these sets are by construction invariant. Within the intersection $\{r_1=\epsilon_1=0\}$, we obtain
\begin{align*}
 v_1' &=-\operatorname{sign}(v_1)
\end{align*}
for any $v_1\ne 0$. Here $v_1=0$ is degenerate (nonsmooth) and this stems from the lack of smoothness along $\mathcal E_{46}$, recall \figref{pp}(a). We resolve this degeneracy by multiplying the right hand side of \eqref{c1eqn} by $\epsilon_1$ and then apply the blowup \eqref{bu1}, using desingularization corresponding to $\bar \epsilon_1$. We will study the dynamics in the chart \eqref{ci2} with $i=1$ in the following section.

\subsection{Analysis in the $\bar \epsilon_1=1$-chart}
Denote the vector-field associated with \eqref{c1eqn} by $V_1$. We then consider the pull-back of $\epsilon_1 V_1$ by the map $(r_1,v_{12},\rho_{12})\mapsto (r_1,v_1,\epsilon_1)$ defined by \eqref{ci2} with $i=1$. This gives the following equations
 \begin{equation}\eqlab{c12eqn}
\begin{aligned}
r_1 '&= {\e^{-(\m-\beta)\epsilon^{-1}}}\left[-1+\e^{(\alpha-\m)\epsilon^{-1}} \e^{-\rho_{12}^{-1}}\e^{v_{12}}+\e^{v_{12}} \right]\e^{-\rho_{12}^{-1}}\rho_{12} r_1,\\
 v_{12}'&={\e^{-v_{12}}-1-\e^{-\rho_{12}^{-1}} },\\
\rho_{12} '&= -{\e^{-(\m-\beta)\epsilon^{-1}}}\left[-1+\e^{(\alpha-\m)\epsilon^{-1}} \e^{-\rho_{12}^{-1}}\e^{v_{12}}+\e^{v_{12}} \right]\e^{-\rho_{12}^{-1}} \rho_{12}^2,
\end{aligned}
 \end{equation}
 after multiplying the right hand side by ${1+\e^{-v_{12}} +\e^{-\rho_{12}^{-1}} }$, which is regular and nonzero within a  fixed compact $D_{12}$ defined by $(r_1,v_{12},\rho_{12})\in [0,r_{10}]\times I\times [0,\rho_{120}]$, with $r_{10},\rho_{120}>0$ both sufficiently small and $I\subset \mathbb R$ a compact interval containing $0$. We take $D_2$ in the $\bar \epsilon=1$-chart large enough so that $D_{12}$ when written in the $(u_2,v_2,r_2)$-coordinates (using the change of coordinates, which is well-defined for $\rho_{12}>0$) intersects $D_2\times [0,\epsilon_0)$ with nonzero interior.  It is clear that $(r_1,v_{12},\rho_{12})=(0,0,0)$ is partially hyperbolic for \eqref{c12eqn} with a single nonzero eigenvalue $-1$. In fact, for $\e^{-(\m-\beta)\epsilon^{-1}}\to 0$ (corresponding to $r_1=0$ or $\rho_{12}=0$) the equation
 \begin{align*}
 \e^{v_{12}} =\frac{1}{1+\e^{-\rho_{12}^{-1}}},\quad \rho_{12}\ge 0,
 \end{align*}
 defines a normally hyperbolic and attracting critical manifold. 

\begin{lemma}\lemmalab{caseceps1lemma}
 The center manifold of $(r_1,v_{12},\rho_{12})=(0,0,0)$ takes the following graph form over $(r_1,\rho_{12})\in [0,r_{10}]\times [0,\rho_{120}]$:
 \begin{align*}
 \e^{v_{12}} &=\frac{1}{1+\e^{-\rho_{12}^{-1}}}+\e^{-\rho_{12}^{-1}}\mathcal O(\e^{-(\m-\beta)\epsilon^{-1}}),
 \end{align*}
 with the remainder term being smooth with respect to $(r_1,\rho_{12})$.
\end{lemma}
\begin{proof}
 This follows from a simple calculation. We leave out further details.
\end{proof}
On the center manifold, we have the following reduced flow
\begin{equation}\nonumber
\begin{aligned}
r_1 '&= {\e^{-(\gamma-\beta)\epsilon^{-1}}}\left[-\frac{1}{1+\e^{-\rho_{12}^{-1} }} +\mathcal O\left(\e^{(\alpha-\gamma)\epsilon^{-1}},\e^{-(\gamma-\beta)\epsilon^{-1} }\right) \right]\e^{-2\rho_{12}^{-1}}\rho_{12} r_1,\\
\rho_{12} '&= -{\e^{-(\gamma-\beta)\epsilon^{-1}}}\left[-\frac{1}{1+\e^{-\rho_{12}^{-1} }} +\mathcal O\left(\e^{(\alpha-\gamma)\epsilon^{-1}},\e^{-(\gamma-\beta)\epsilon^{-1} }\right) \right]\e^{-2\rho_{12}^{-1}} \rho_{12}^2,
\end{aligned}
 \end{equation}
 Here we have used \eqref{mepsdef2}.
It follows that $r_1$ is decreasing for $r_1,\rho_{12}>0$.

\begin{remark}\remlab{caseceps1extrem}
In the chart $\bar \epsilon_1$-chart, we have $\epsilon = r_1\rho_{12}$. For $r_{1}\ge r_{10}>0$, we can therefore replace $\rho_{12}$ by $\rho_{12}=r_1^{-1}\epsilon$. This leads to the following system
 \begin{equation}\nonumber
\begin{aligned}
r_1 '&= {\e^{-(\m-\beta)\epsilon^{-1}}}\left[-1+\e^{(\alpha-\m)\epsilon^{-1}} \e^{-\epsilon^{-1}r_1}\e^{v_{12}}+\e^{v_{12}} \right]\e^{-\epsilon^{-1} r_1 }\epsilon,\\
 v_{12}'&={\e^{-v_{12}}-1-\e^{-\epsilon^{-1} r_1}},\\
 \epsilon ' &=0,
\end{aligned}
 \end{equation}
 which is slow-fast. This system can also be obtained by blowing up $\mathcal E_{46}\times\{0\}$ directly for \eqref{uvext} (away from $u=\gamma$) using
 \begin{align*}
  r\ge 0,\,(\bar v,\bar \epsilon)\in S^1\mapsto \begin{cases}
                                                 v = r\bar v,\\
                                                 \epsilon =r\bar \epsilon,
                                                \end{cases}
 \end{align*}
and the associated $\bar \epsilon=1$-chart.

 From here we also see that Fenichel's theory applies uniformly with respect to $r_1\ge r_{10}>0$. Indeed, we can compactify $r_1\ge r_{10}$ by considering the system in terms of $(r_1^{-1},v_{12},\epsilon)$. We then see that we can take the unique manifold from infinity,  recall the discussion around \eqref{zv} in the introduction, as a slow manifold for all $0<\epsilon\ll 1$. The statement in \lemmaref{caseceps1lemma} therefore provides a continuation of the invariant manifold from infinity (as a slow manifold) through the invariant foliations $\epsilon=r_1\rho_{12}=\text{const}$.

\end{remark}
\subsection{Analysis in the $\bar u=-1$-chart}
Next, in the $\bar u=-1$-chart, we fix a compact set $D_3\subset \mathbb R^3$ defined by $(r_3,v_3,\epsilon_3)\in [0,r_{30})\times I\times [0,\epsilon_{30}]$, with $r_{30},\epsilon_{30}>0$ both sufficiently small and $I\subset \mathbb R$ a compact interval containing $0$. The pull-back of the vector-field \eqref{uvext} by the map $(r_3,v_3,\epsilon_3)\mapsto (u,v,\epsilon)$ defined by \eqref{c3}, then takes the following form
\begin{equation}\eqlab{c3eqn}
\begin{aligned}
r_3 '&= -\frac{\e^{-(\m-\beta)\epsilon^{-1}}}{1+\e^{-\epsilon_3^{-1}(1+v_3)}}\left[-1+\e^{(\alpha-\m)\epsilon^{-1}} \e^{\epsilon_3^{-1}}\e^{\epsilon_3^{-1} v_3}+\e^{\epsilon_3^{-1} v_3} \right]r_3,\\
 v_3'&=\frac{\e^{-\epsilon_3^{-1}(1+v_3)}-1-\e^{-\epsilon_3^{-1}} }{1+\e^{-\epsilon_3^{-1}(1+v_3)}}\\
 &+\frac{\e^{-(\m-\beta)\epsilon^{-1}}}{1+\e^{-\epsilon_3^{-1}(1+v_3)}}\left[-1+\e^{(\alpha-\m)\epsilon^{-1}} \e^{\epsilon_3^{-1}}\e^{\epsilon_3^{-1} v_3}+\e^{\epsilon_3^{-1} v_3} \right]v_3,\\
\epsilon_3 '&= \frac{\e^{-(\m-\beta)\epsilon^{-1}}}{1+\e^{-\epsilon_3^{-1}(1+v_3)}}\left[-1+\e^{(\alpha-\m)\epsilon^{-1}} \e^{\epsilon_3^{-1}}\e^{\epsilon_3^{-1} v_3}+\e^{\epsilon_3^{-1} v_3} \right]\epsilon_3.
\end{aligned}
 \end{equation}
Here we have as discussed above used a desingularization corresponding to division of the right hand side by $\epsilon_3$. Moreover, we have for convenience also multiplied the right hand side by
\begin{align*}
 \frac{\e^{-(\m-\beta)\epsilon^{-1}}\sum_{i=1}^6 e^{F_i(\m-r_3,r_3 v_3)\epsilon^{-1}} }{1+\e^{-\epsilon_3^{-1}(1+v_3)}},
\end{align*}
which is regular and nonzero on $D_3$. As above, we have left $\epsilon=r_3\epsilon_3$ as a constant where it makes sense. The system clearly extends smoothly to $\{r_3=0\}$ and $\{\epsilon_3=0\}$ provided that $v_3\ne -1$. In fact, these sets are by construction invariant. Within $\{r_3=0\}$ we have
\begin{equation}\nonumber
\begin{aligned}
 v_3'&=\frac{\e^{-\epsilon_3^{-1}(1+v_3)}-1-\e^{-\epsilon_3^{-1}} }{1+\e^{-\epsilon_3^{-1}(1+v_3)}},\\
\epsilon_3 '&= 0,
\end{aligned}
 \end{equation}
and therefore within the intersection $\{r_3=\epsilon_3=0\}$, we obtain that
%
%
\begin{align*}
 v_3' &=-\operatorname{sign}(1+v_3)
\end{align*}
for any $v_3\ne -1$. Here $v_3=-1$ is degenerate (nonsmooth) and this stems from the lack of smoothness along $\mathcal E_{45}$, recall \figref{pp}(a). We resolve this degeneracy by multiplying the right hand side of \eqref{c3eqn} by $\epsilon_3$ and then apply the blowup \eqref{bu3}, using desingularization corresponding to $\bar \epsilon_3$. We will study the dynamics in the chart \eqref{ci2} with $i=3$ in the following section.
\subsection{Analysis in the $\bar \epsilon_3=1$-chart}
Denote the vector-field associated with \eqref{c3eqn} by $V_3$. We then consider the pull-back of $\epsilon_3 V_3$ by the map $(r_3,v_{32},\rho_{32})\mapsto (r_3,v_3,\epsilon_3)$ defined by \eqref{ci2} with $i=3$. This gives the following equations
 \begin{equation}\eqlab{c32eqn}
\begin{aligned}
r_3 '&= -\e^{-(\m-\beta)\epsilon^{-1}}\left[-1+\e^{(\alpha-\m)\epsilon^{-1}} \e^{v_{32}}+\e^{-\rho_{32}^{-1}}\e^{v_{32}} \right]r_3\rho_{32},\\
 v_{32}'&={\e^{-v_{32}}-1-\e^{-\rho_{32}^{-1}} }\\
 &-\e^{-(\m-\beta)\epsilon^{-1}}\left[-1+\e^{(\alpha-\m)\epsilon^{-1}} \e^{v_{32}}+\e^{-\rho_{32}^{-1}}\e^{v_{32}} \right],\\
\rho_{32} '&= \e^{-(\m-\beta)\epsilon^{-1}}\left[-1+\e^{(\alpha-\m)\epsilon^{-1}} \e^{v_{32}}+\e^{-\rho_{32}^{-1}}\e^{v_{32}} \right]\rho_{32}^2.
\end{aligned}
 \end{equation}
 after multiplying the right hand side by $1+\e^{-v_{32}}$, which is regular and nonzero within a  fixed compact $D_{32}$ defined by $(r_3,v_{32},\rho_{32})\in [0,r_{30}]\times I\times [0,\rho_{320}]$, with $r_{30},\rho_{320}>0$ both sufficiently small and $I\subset \mathbb R$ a compact interval containing $0$.  We take $D_2$ in the $\bar \epsilon=1$-chart large enough so that $D_{32}$ when written in the $(u_2,v_2,r_2)$-coordinates (using the change of coordinates, which is well-defined for $\rho_{32}>0$) intersects $D_2\times [0,\epsilon_0)$ with nonzero interior. It is clear that $(r_3,v_{32},\rho_{32})=(0,0,0)$ is partially hyperbolic for \eqref{c32eqn} with a single nonzero eigenvalue $-1$. In fact, for $\e^{-(\m-\beta)\epsilon^{-1}}\to 0$ (corresponding to $r_3=0$ or $\rho_{32}=0$) the equation
 \begin{align*}
 \e^{v_{32}} =\frac{1}{1+\e^{-\rho_{32}^{-1}}},\quad \rho_{32}\ge 0,
 \end{align*}
  defines a normally hyperbolic and attracting critical manifold.

\begin{lemma}\lemmalab{caseceps3lemma}
 The center manifold of $(r_3,v_{32},\rho_{32})=(0,0,0)$ takes the following graph form over $(r_3,\rho_{32})\in [0,r_{30}]\times [0,\rho_{320}]$:
 \begin{align*}
 \e^{v_{32}} &=\frac{1}{1+\e^{-\rho_{32}^{-1}}}+\mathcal O(\e^{-(\m-\beta)\epsilon^{-1}}),
 \end{align*}
 with the remainder term being smooth with respect to $(r_3,\rho_{32})$.
\end{lemma}
\begin{proof}
 This follows from a simple calculation. We leave out further details.
\end{proof}
On the center manifold, we have the following reduced flow
\begin{equation}\nonumber
\begin{aligned}
r_3 '&= -\e^{-(\gamma-\beta)\epsilon^{-1}}\left[-1+\mathcal O\left(\e^{-\rho_{32}^{-1}},\e^{(\alpha-\gamma)\epsilon^{-1}}\right)\right]r_3\rho_{32},\\
\rho_{32} '&= \e^{-(\gamma-\beta)\epsilon^{-1}}\left[-1+\mathcal O\left(\e^{-\rho_{32}^{-1}},\e^{(\alpha-\gamma)\epsilon^{-1}}\right)\right]\rho_{32}^2.
\end{aligned}
 \end{equation}
 Here we have again used \eqref{mepsdef2}.
It follows that $r_3$ is increasing for $r_3,\rho_{32}>0$.


\begin{remark}\remlab{caseceps3extrem}
In the chart $\bar \epsilon_3=1$, we have $\epsilon = r_3\rho_{32}$. For $r_{3}\ge r_{30}>0$, we can therefore replace $\rho_{32}$ by $\rho_{32}=r_3^{-1}\epsilon$. This leads to the following system
\begin{equation}\nonumber
\begin{aligned}
r_3 '&= -\e^{-(\m-\beta)\epsilon^{-1}}\left[-1+\e^{(\alpha-\m)\epsilon^{-1}} \e^{v_{32}}+\e^{-\epsilon^{-1}r_3}\e^{v_{32}} \right]\epsilon,\\
 v_{32}'&={\e^{-v_{32}}-1-\e^{-\epsilon^{-1}r_3  }}\\
 &-\e^{-(\m-\beta)\epsilon^{-1}}\left[-1+\e^{(\alpha-\m)\epsilon^{-1}} \e^{v_{32}}+\e^{-\epsilon^{-1}r_3}\e^{v_{32}} \right],\\
 \epsilon' &=0,
\end{aligned}
 \end{equation}
  which is slow-fast. This system can also be obtained by blowing up $\mathcal E_{45}\times \{0\}$ directly for \eqref{uvext} (away from $u=\gamma$) using
  \begin{align*}
   r\ge 0,\,(\bar v,\bar \epsilon)\in S^1\mapsto \begin{cases}
                                                  v = u-\m+r\bar v,\\
                                                  \epsilon =r\bar \epsilon,
                                                 \end{cases}
  \end{align*}
and the associated $\bar \epsilon=1$-chart.

From here we see that the invariant manifold in \lemmaref{caseceps3lemma} can be extended as slow manifolds uniformly with respect to $r_3\ge r_{30}>0$ for all $0<\epsilon\ll 1$. Indeed, we can compactify $r_3\ge r_{30}$ by considering the system in terms of $(r_3^{-1},v_{32},\epsilon)$. Notice that $r_3\rightarrow \infty$ corresponds to $u\to -\infty$ (which corresponds to $\tilde u\to 0$, see \eqref{sc}).

\end{remark}

In conclusion, we have shown that there is a normally hyperbolic and attracting critical manifold on the blown up space (in a neighborhood of $(u,v,\epsilon)=(\gamma,0,0)$). We summarize the results of the section in \figref{ppcbulocal2}.

 \begin{figure}
  \begin{center}
    {\includegraphics[width=0.475\textwidth]{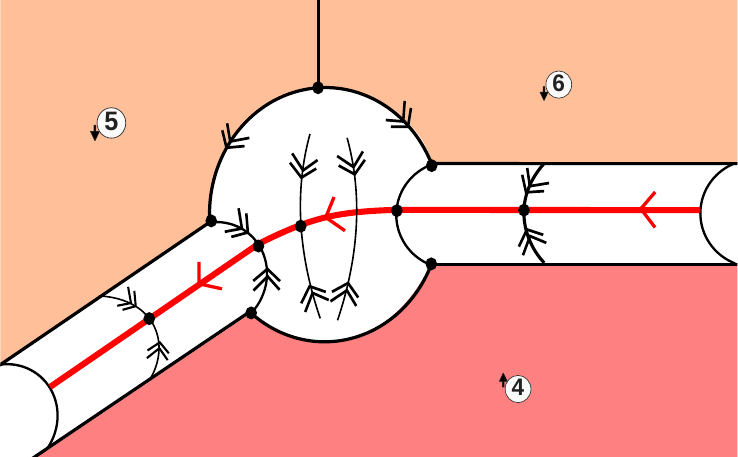}}
%
    \caption{Illustration of the results in the case (i) on the blowup of $(u,v,\epsilon)=(c,0,0)$. The normally hyperbolic and attracting critical manifold is in red. Notice that the third dimension coming out of the page corresponds to $\epsilon$. }\figlab{ppcbulocal2}
        \end{center}
    \end{figure}
\section{Case (iii)}\seclab{casea}
In this section, we consider case (iii): $\max(\alpha,\gamma) <\beta$, $\alpha\ne \gamma$.
The corresponding tropical dynamics is illustrated in \figref{pp}(c). We again consider \eqref{uvext}
and now blowup $(u,v,\epsilon)=(\beta,0,0)$ using the following spherical blowup transformation:
\begin{align}\eqlab{bucasea1}
 r\ge 0,\,(\bar u,\bar v,\bar \epsilon)\in S^2\mapsto \begin{cases}
                                                       u =\beta+r\bar u,\\
                                                       v =r\bar v,\\
                                                       \epsilon =r\bar \epsilon.
                                                      \end{cases}
\end{align}
We work with the following charts (redefining the chart-specific coordinates):
\begin{align}
 \bar u&=1:\quad \begin{cases}
                 u =\beta+r_1,\\
                 v =r_1 v_1,\\
                 \epsilon =r_1\epsilon_1,
                \end{cases}\eqlab{c1casea}\\
                \bar \epsilon&=1:\quad \begin{cases}
                 u =\beta+r_2u_2,\\
                 v =r_2 v_2,\\
                 \epsilon =r_2,
                \end{cases}\eqlab{c2casea}\\
                \bar u&=-1:\quad \begin{cases}
                 u =\beta-r_3,\\
                 v =r_3 v_3,\\
                 \epsilon =r_3\epsilon_3.
                \end{cases}\eqlab{c3casea}
\end{align}
As above, we achieve the desingularization upon dividing the pull-back vector-field by $\bar \epsilon$; we use $ \epsilon_1$, $1$ resp. $ \epsilon_3$ in the three charts.

\begin{remark}
 Notice that the $\bar \epsilon=1$-chart corresponds to
 \begin{align*}
 \begin{cases}
  \tilde u = \tilde \beta \e^{u_2},\\
\tilde v = \e^{v_2}.
\end{cases}
 \end{align*}
 with respect to the original parameters and variables (recall \eqref{abc} and \eqref{sc}).

\end{remark}

The change of coordinates are defined by the expressions
\begin{align*}
 \begin{cases}
  r_1 =r_2 u_2,\\
  v_1 =u_2^{-1} v_2,\\
  \epsilon_1 =u_2^{-1},
 \end{cases}\quad
\begin{cases}
  r_3 =-r_2 u_2,\\
  v_3 =(-u_2)^{-1} v_2,\\
  \epsilon_3 =(-u_2)^{-1},
 \end{cases}
\end{align*}
for $u_2>0$ resp. $u_2<0$.

In the charts $\bar u=\pm 1$, we find additional degeneracies due to the nonsmoothness along $\mathcal E_{13}$ and $\mathcal E_{46}$, see \figref{pp}(c), that we (as above, following \cite{jelbart2025switchingmultipletimescalesgeometric,kaklamanos2019a}) resolve using additional cylindrical blowups. For example, in $\bar u=1$, we blow up the degenerate line $r_1\ge 0$, $(v_1,\epsilon_1)=(0,0)$, as follows
 \begin{align}\eqlab{bu1casea}
 \rho \ge 0,\,(\bar v_1,\bar \epsilon_1)\in S^1\mapsto \begin{cases}
                                                       v_1 =\rho \bar v_1,\\
                                                       \epsilon_1 =\rho \bar \epsilon_1,
                                                      \end{cases}
\end{align}
leaving $r_1$ untouched.
Similarly, in $\bar u=-1$, we blow up the line $r_3\ge 0$, $(v_3,\epsilon_3)=(0,0)$ as follows:
\begin{align}\eqlab{bu3casea}
 \rho \ge 0,\,(\bar v_3,\bar \epsilon_3)\in S^1\mapsto \begin{cases}
                                                       v_3 =\rho \bar v_3,\\
                                                       \epsilon_3 =\rho \bar \epsilon_3,
                                                      \end{cases}
\end{align}
leaving $r_3$ untouched.
For simplicity, we will restrict attention to the following charts
\begin{equation}\eqlab{ci2casea}
\begin{aligned}
                 \bar \epsilon_i&=1:\quad \begin{cases}
                 v_i =\rho_{i2}v_{i2},\\
                 \epsilon_i =\rho_{i2},
                \end{cases}
\end{aligned}
\end{equation}
for $i\in \{1,3\}$, resp. Notice that upon composing this with \eqref{c1casea} resp. \eqref{c3casea}, we find that
\begin{align}\eqlab{vv32}
v  &= \epsilon v_{i2}.
\end{align}
Hence $v_2=v_{i2}$, see \eqref{c2casea}.





We illustrate the result of the blowup sequence in \figref{ppabulocal}. Notice in comparison with \figref{ppau}, that our analysis here only limits to the sets relevant to slow manifolds.  The details of resolving the remaining degeneracies are similar and therefore left out for simplicity.

 \begin{figure}
  \begin{center}
    {\includegraphics[width=0.875\textwidth]{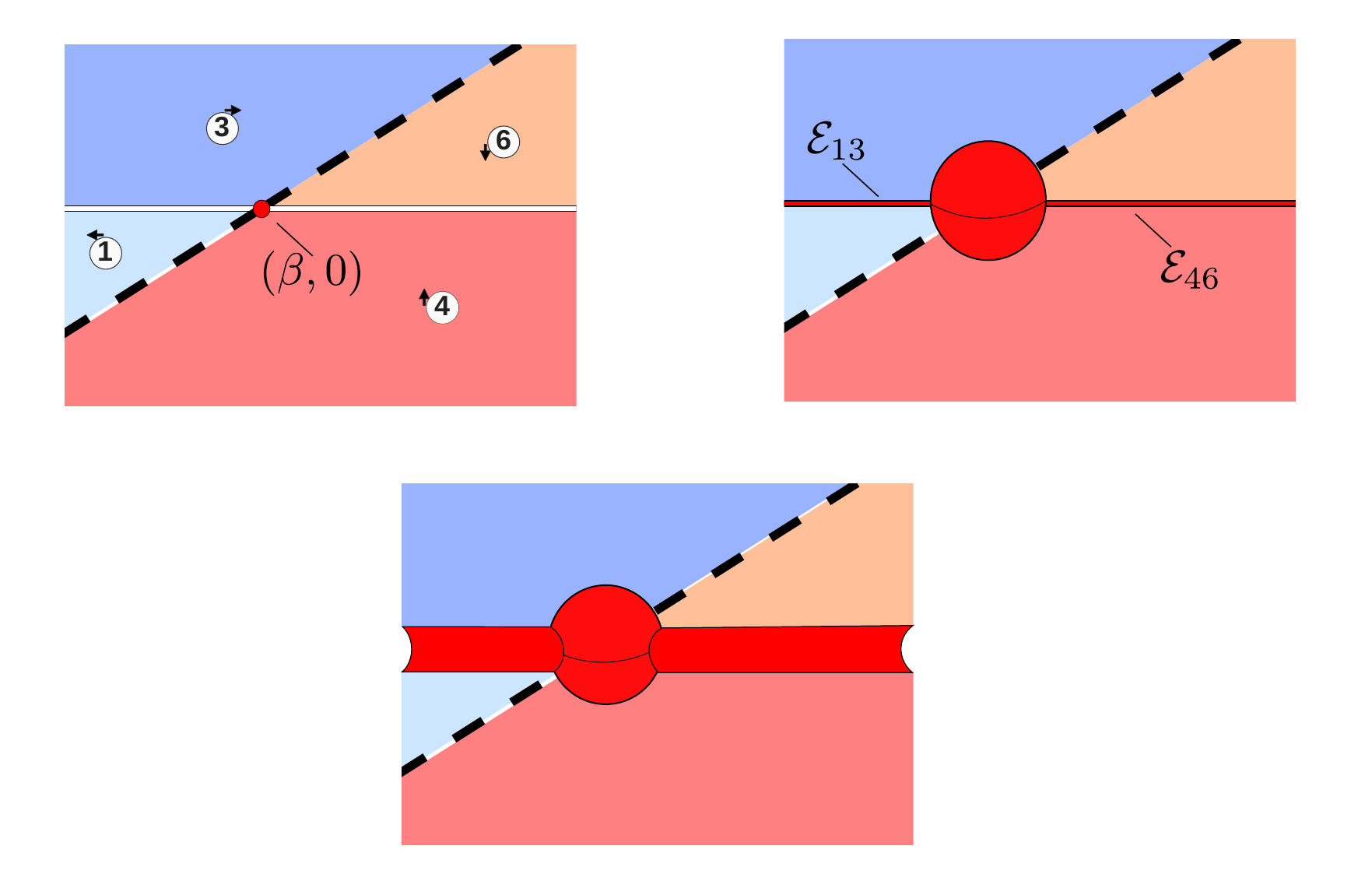}}
%
    \caption{Illustration of the blowup sequence defined by \eqref{bucasea1} and \eqref{bu1casea}, \eqref{bu3casea}, going clockwise starting from the top left. First we blow up $(u,v,\epsilon)=(\beta,0,0)$, indicated in the top left by a red dot in the $(u,v)$-plane, to a sphere (top right). Notice that the third dimension coming out of the page corresponds to $\epsilon$. Subsequently, we blow up the lines $\mathcal E_{13} \times \{0\}$ and $\mathcal E_{46}\times \{0\}$ to cylinders (in red). The details of the blowups of the dotted lines (corresponding to $\mathcal E_{14}\times \{0\}$ and $\mathcal E_{36}\times \{0\}$) to separate cylinders are left out for simplicity as this does not relate to slow manifolds. }\figlab{ppabulocal}
        \end{center}
    \end{figure}


\subsection{Analysis in the $\bar \epsilon=1$-chart}\seclab{bareps1}
We first consider the $\bar \epsilon=1$-chart. We fix a sufficiently large compact domain $D_2\subset \mathbb R^2$ and consider the pull-back of the vector-field \eqref{uvext} by the map $(u_2,v_2,r_2)\mapsto (u,v,\epsilon)$ defined by \eqref{c2casea}. This leads to the following equations
\begin{equation}\eqlab{uv2casea}
\begin{aligned}
 u_2' &=\e^{v_2}-1+\e^{(\alpha-\beta)\epsilon^{-1}}\e^{v_2-u_2},\\
 v_2' &= \e^{u_2}(\e^{-v_2}-1)-\e^{(\m-\beta)\epsilon^{-1}},
\end{aligned}
\end{equation}
and $r_2'=\epsilon'=0$, after multiplication of the right hand side by $$\sum_{i=1}^6 \e^{F_i(\beta +\epsilon u_2,\epsilon v_2)\epsilon^{-1}}\to 1+\e^{u_2}+\e^{v_2}+\e^{u_2-v_2}\quad \mbox{for}\quad \epsilon\to 0,$$ which is regular and nonzero for $(u_2,v_2)\in D_2$ as $\epsilon\rightarrow 0$. (The time used in \eqref{uv2casea} is therefore also the original one in \eqref{mm}.) For $\epsilon\to 0$, we obtain
\begin{equation}\nonumber
\begin{aligned}
 u_2' &=\e^{v_2}-1,\\
 v_2' &= \e^{u_2}(\e^{-v_2}-1),
\end{aligned}
\end{equation}
since $\max(\alpha,\gamma)<\beta$, recall also \eqref{meps}.
\begin{lemma}\lemmalab{bareps1lemma}
 The system \eqref{uv2casea} is slow-fast with respect to $\epsilon\to 0$ with a normally hyperbolic and attracting critical manifold $C_2$ defined by
 \begin{align*}
  v_2=0.
 \end{align*}
We have the following expansion of an associated slow manifold $S_{2,\epsilon}$ with $S_{2,0}=C_2\cap D_2$:
\begin{align*}
  \e^{v_2} = \frac{\e^{u_2}}{\e^{u_2}+\e^{(\alpha-\beta)\epsilon^{-1}}}-\e^{-u_2}\e^{(\gamma-\beta)\epsilon^{-1}} \left(1+\mathcal O\left(\e^{(\alpha-\beta)\epsilon^{-1}},\e^{(\gamma-\beta)\epsilon^{-1}}\right)\right),
\end{align*}
for all $0\le \epsilon<\epsilon_0$, with $\epsilon_0>0$ sufficiently small, and where the remainder is smooth with respect $(\e^{u_2},\e^{(\alpha-\beta)\epsilon^{-1}},\e^{(\gamma-\beta)\epsilon^{-1}})$.
\end{lemma}
\begin{proof}
 The first statement is obvious. Now for the expansion of the slow-manifold, we define $\tilde u_2=\e^{u_2}$, $\tilde v_2=\e^{v_2}$ so that
   \begin{equation}\eqlab{x2y2casea}
\begin{aligned}
 \tilde u_2' &=-\left(\tilde u_2(1-\tilde v_2)-\e^{(\alpha-\beta)\epsilon^{-1}}\tilde v_2\right),\\
 \tilde v_2' &= \tilde u_2(1-\tilde v_2)-\e^{(\alpha-\beta)\epsilon^{-1}}\tilde v_2-\e^{(\gamma-\beta)\epsilon^{-1}}\tilde v_2,
\end{aligned}
\end{equation}
Here we have also used the definition of $\m$, recall \eqref{mepsdef}. Let $\e^{(\gamma-\beta)\epsilon^{-1}}\to 0$. Then
\begin{align*}
\tilde v_2 &=\frac{\tilde u_2}{\tilde u_2+\e^{(\alpha-\beta)\epsilon^{-1}}} = 1+\mathcal O(\e^{(\alpha-\beta)\epsilon^{-1}}),
\end{align*}
with $\tilde u_2$ uniformly bounded away from $0$,
is a normally hyperbolic and attracting critical manifold. It follows that the slow manifold can be written as an asymptotic expansion in $\e^{(\gamma-\beta)\epsilon^{-1}}$ with $(\tilde u_2,\e^{(\alpha-\beta)\epsilon^{-1}})$-dependent coefficients.
The result then follows from a simple calculation.
\end{proof}

On $S_{2,\epsilon}$ we have the following slow flow
\begin{align*}
 u_2'&=-\e^{-u_2} \e^{(\gamma-\beta)\epsilon^{-1}} \left(1+\mathcal O\left(\e^{(\alpha-\beta)\epsilon^{-1}},\e^{(\gamma-\beta)\epsilon^{-1}}\right)\right),
\end{align*}
with the remainder being smooth with respect $(\e^{u_2},\e^{(\alpha-\beta)\epsilon^{-1}},\e^{(\gamma-\beta)\epsilon^{-1}})$. Hence $u_2$ is decreasing for all $0<\epsilon\ll 1$.
\subsection{Analysis in the $\bar u=1$-chart}
Next, in the $\bar u=1$-chart, we fix a compact set $D_1\subset \mathbb R^3$ defined by $(r_1,v_1,\epsilon_1)\in [0,r_{10})\times I\times [0,\epsilon_{10}]$, with $r_{10},\epsilon_{10}>0$ both sufficiently small and $I\subset \mathbb R$ a compact interval containing $0$. The pull-back of the vector-field \eqref{uvext} by the map $(r_1,v_1,\epsilon_1)\mapsto (u,v,\epsilon)$ defined by \eqref{c1casea}, then takes the following form
\begin{equation}\eqlab{c1eqncasea}
\begin{aligned}
r_1 '&= \frac{\e^{-\epsilon_1^{-1} }}{1+\e^{\epsilon_1^{-1}(v_1-1)}+ \e^{-\epsilon_1^{-1} v_1} }\left[-1+\e^{(\alpha-\beta)\epsilon^{-1}} \e^{\epsilon_1^{-1}(v_1-1) }+\e^{\epsilon_1^{-1} v_1} \right]r_1,\\
 v_1'&=\frac{\e^{-\epsilon_1^{-1} v_1} -1-\e^{(\m-\beta)\epsilon^{-1}}\e^{-\epsilon_1^{-1}}  }{1+\e^{\epsilon_1^{-1}(v_1-1)}+ \e^{-\epsilon_1^{-1} v_1} }\\
 &-\frac{\e^{-\epsilon_1^{-1} }}{1+\e^{\epsilon_1^{-1}(v_1-1)}+ \e^{-\epsilon_1^{-1} v_1} }\left[-1+\e^{(\alpha-\beta)\epsilon^{-1}} \e^{\epsilon_1^{-1}(v_1-1) }+\e^{\epsilon_1^{-1} v_1} \right]v_1,\\
\epsilon_1 '&= -\frac{\e^{-\epsilon_1^{-1} }}{1+\e^{\epsilon_1^{-1}(v_1-1)}+ \e^{-\epsilon_1^{-1} v_1} }\left[-1+\e^{(\alpha-\beta)\epsilon^{-1}} \e^{\epsilon_1^{-1}(v_1-1) }+\e^{\epsilon_1^{-1} v_1} \right]\epsilon_1.
\end{aligned}
 \end{equation}
Here we have as discussed above used a desingularization corresponding to division of the right hand side by $\epsilon_1$. Moreover, we have for convenience also multiplied the right hand side by
\begin{align*}
 \frac{\e^{-\epsilon_1^{-1}} \sum_{i=1}^6 e^{F_i(\beta+r_1,r_1 v_1)\epsilon^{-1}} }{1+\e^{\epsilon_1^{-1}(v_1-1)}+ \e^{-\epsilon_1^{-1} v_1} }\to 1\quad \mbox{for}\quad (r_1,\epsilon_1)\rightarrow (0,0),
\end{align*}
which is regular and nonzero on $D_1$. We have left $\epsilon=r_1\epsilon_1$ as a constant where it makes sense. The system clearly extends smoothly to $\{r_1=0\}$ and $\{\epsilon_1=0\}$ provided that $v_1\ne 1$ and $v_1\ne 0$. In fact, these sets are by construction invariant. Within $\{r_1=0\}$ we have
\begin{equation}\nonumber
\begin{aligned}
 v_1'&=\frac{\e^{-\epsilon_1^{-1} v_1} -1  }{1+\e^{\epsilon_1^{-1}(v_1-1)}+ \e^{-\epsilon_1^{-1} v_1} }\\
 &-\frac{\e^{-\epsilon_1^{-1} }}{1+\e^{\epsilon_1^{-1}(v_1-1)}+ \e^{-\epsilon_1^{-1} v_1} }\left[-1+\e^{\epsilon_1^{-1} v_1} \right]v_1,\\
\epsilon_1 '&= -\frac{\e^{-\epsilon_1^{-1} }}{1+\e^{\epsilon_1^{-1}(v_1-1)}+ \e^{-\epsilon_1^{-1} v_1} }\left[-1+\e^{\epsilon_1^{-1} v_1} \right]\epsilon_1,
\end{aligned}
 \end{equation}
and therefore within $\{r_1=\epsilon_1=0\}$, we obtain
\begin{align*}
 v_1' &=-\operatorname{sign}(v_1)
\end{align*}
for any $v_1\ne 1$ and $v_1\ne 0$. The degeneracies along $v_1=0$ and $v_1=1$ stem from the lack of smoothness along $\mathcal E_{46}$ and $\mathcal E_{36}$, see \figref{pp}(c). We resolve the degeneracy along $\mathcal E_{46}$ by multiplying the right hand side of \eqref{c1eqncasea} by $\epsilon_1$ and then apply the blowup \eqref{bu1casea}, using desingularization corresponding to $\bar \epsilon_1$. We will study the dynamics in the chart \eqref{ci2casea} with $i=1$ in the following section. For simplicity we leave out the details of $\mathcal E_{36}$ as it does not relate to slow manifolds.

\subsection{Analysis in the $\bar \epsilon_1=1$-chart}\seclab{bareps11}
Denote the vector-field associated with \eqref{c1eqncasea} by $V_1$. We then consider the pull-back of $\epsilon_1 V_1$ by the map $(r_1,v_{12},\rho_{12})\mapsto (r_1,v_1,\epsilon_1)$ defined by \eqref{ci2casea} with $i=1$. This gives the following equations
 \begin{equation}\eqlab{c12eqncasea}
\begin{aligned}
r_1 '&= \e^{-\rho_{12}^{-1} }\left[-1+\e^{(\alpha-\beta)\epsilon^{-1}} \e^{-\rho_{12}^{-1}}\e^{v_{12}} +\e^{v_{12}} \right]\rho_{12} r_1,\\
 v_{12}'&=\e^{-v_{12}}-1-\e^{(\m-\beta)\epsilon^{-1}}\e^{-\rho_{12}^{-1}},\\
\rho_{12} '&= -\e^{-\rho_{12}^{-1} }\left[-1+\e^{(\alpha-\beta)\epsilon^{-1}} \e^{-\rho_{12}^{-1}}\e^{v_{12}} +\e^{v_{12}} \right]\rho_{12}^2,
\end{aligned}
 \end{equation}
 after multiplying the right hand side by $1+\e^{-\rho_{12}^{-1}}\e^{v_{12}}+ \e^{-v_{12}}$, which is regular and nonzero within a  fixed compact $D_{12}$ defined by $(r_1,v_{12},\rho_{12})\in [0,r_{10}]\times I\times [0,\rho_{120}]$, with $r_{10},\rho_{120}>0$ both sufficiently small and $I\subset \mathbb R$ a compact interval containing $0$. We take $D_2$ in the $\bar \epsilon=1$-chart large enough so that $D_{12}$ when written in the $(u_2,v_2,r_2)$-coordinates (using the change of coordinates, which is well-defined for $\rho_{12}>0$) intersects $D_2\times [0,\epsilon_0)$ with nonzero interior.

 It is clear that $v_{12}=0$ defines a partially hyperbolic critical manifold within either $r_1=0$ or $\rho_{12}=0$. In particular, the linearization of $(r_1,v_{12},\rho_{12})=(0,0,0)$ has a single nonzero eigenvalue $-1$.

%

\begin{lemma}\lemmalab{caseaeps1lemma}
 The center manifold of $(r_1,v_{12},\rho_{12})=(0,0,0)$ takes the following graph form over $(r_1,\rho_{12})\in [0,r_{10}]\times [0,\rho_{120}]$:
 \begin{align*}
 \e^{v_{12}} &=\frac{1}{1+\e^{(\alpha-\beta)\epsilon^{-1}} \e^{-\rho_{12}^{-1}}}-\e^{-\rho_{12}^{-1}}\e^{(\gamma-\beta)\epsilon^{-1}}\left(1+\mathcal O\left(\e^{(\alpha-\beta)\epsilon^{-1}},\e^{(\gamma-\beta)\epsilon^{-1}}\right)\right),
 \end{align*}
 with the remainder term being smooth with respect to $(r_1,\rho_{12})$.
\end{lemma}
\begin{proof}
 We first use that for $\e^{(\gamma-\beta)\epsilon^{-1}}\to 0$ then
 \begin{align}
  \e^{v_{12}} &=\frac{1}{1+\e^{(\alpha-\beta)\epsilon^{-1}} \e^{-\rho_{12}^{-1}}},\eqlab{cm12casea0}
 \end{align}
defines a critical manifold.  Similarly, by augmenting an equation for $q=\e^{-\rho_{12}^{-1}}$:
\begin{align}\eqlab{qq}
 q' = -q^2 \left[-1+\e^{(\alpha-\beta)\epsilon^{-1}} \e^{-\rho_{12}^{-1}}\e^{v_{12}} +\e^{v_{12}} \right],
\end{align}
using $q'=\epsilon^{-1} \rho_{12}^{-2} q \rho_{12}'$,
we see that \eqref{cm12casea0} defines a critical manifold for the extended system with respect to $(r_1,v_{12},\rho_{12},q)$, obtained from \eqref{c12eqncasea} (with all $\e^{-\rho_{12}^{-1}}$ replaced by $q$) and \eqref{qq}, within $q=0$. Consequently, the center manifold is an asymptotic expansion in $\e^{(\gamma-\beta)\epsilon^{-1}}$ and $q=\e^{-\rho_{12}^{-1}}$ with coefficients that depend upon $(r_1,\rho_{12},\e^{(\alpha-\beta)\epsilon^{-1}})$. A simple calculation then gives the result.
\end{proof}
On the center manifold, we have the following reduced flow
\begin{equation}\nonumber
\begin{aligned}
r_1 '&= {\e^{-(\gamma-\beta)\epsilon^{-1}}}\e^{-2\rho_{12}^{-1}}\left[-1 +\mathcal O\left(\e^{(\alpha-\beta)\epsilon^{-1}},\e^{-(\gamma-\beta)\epsilon^{-1} }\right) \right]\rho_{12} r_1,\\
\rho_{12} '&=- {\e^{-(\gamma-\beta)\epsilon^{-1}}}\e^{-2\rho_{12}^{-1}}\left[-1 +\mathcal O\left(\e^{(\alpha-\beta)\epsilon^{-1}},\e^{-(\gamma-\beta)\epsilon^{-1} }\right) \right]\rho_{12}^2,
\end{aligned}
 \end{equation}
It follows that $r_1$ is decreasing for $r_1,\rho_{12}>0$.

\begin{remark}\remlab{caseaeps1extrem}
In the chart $\bar \epsilon_1=1$, we have $\epsilon = r_1\rho_{12}$. For $r_{1}\ge r_{10}>0$, we can therefore replace $\rho_{12}$ by $\rho_{12}=r_1^{-1}\epsilon$. This leads to the following system
\begin{equation}\nonumber
\begin{aligned}
r_1 '&= \e^{-\epsilon^{-1}r_1 }\left[-1+\e^{(\alpha-\beta)\epsilon^{-1}} \e^{-\epsilon^{-1}r_1}\e^{v_{12}} +\e^{v_{12}} \right]\epsilon,\\
 v_{12}'&=\e^{-v_{12}}-1-\e^{(\m-\beta)\epsilon^{-1}}\e^{-\epsilon^{-1}r_1},\\
 \epsilon'&=0,
\end{aligned}
 \end{equation}
 which is slow-fast. This system can also be obtained by blowing up $\mathcal E_{46}\times\{0\}$ directly for \eqref{uvext} (away from $u=\beta$) using
 \begin{align*}
  r\ge 0,\,(\bar v,\bar \epsilon)\in S^1\mapsto \begin{cases}v =r\bar v,\\
                                                 \epsilon =r\bar \epsilon,
                                                \end{cases}
 \end{align*}
and the associated $\bar \epsilon=1$-chart.

From here we also see that Fenichel's theory applies uniformly with respect to $r_1\ge r_{10}>0$. Indeed, we can compactify $r_1\ge r_{10}$ by considering the system in terms of $(r_1^{-1},v_{12},\epsilon)$. We then see that we can take the unique manifold from infinity,  recall the discussion around \eqref{zv} in the introduction, as a slow manifold for all $0<\epsilon\ll 1$. The statement \lemmaref{caseaeps1lemma} therefore provides a continuation of the invariant manifold from infinity (as a slow manifold) through the invariant foliations $\epsilon=r_1\rho_{12}=\textrm{const}$.

\end{remark}

\subsection{Analysis in the $\bar u=-1$-chart}
Next, in the $\bar u=-1$-chart, we fix a compact set $D_3\subset \mathbb R^3$ defined by $(r_3,v_3,\epsilon_3)\in [0,r_{30})\times I\times [0,\epsilon_{30}]$, with $r_{30},\epsilon_{30}>0$ both sufficiently small and $I\subset \mathbb R$ a compact interval containing $0$. The pull-back of the vector-field \eqref{uvext} by the map $(r_3,v_3,\epsilon_3)\mapsto (u,v,\epsilon)$ defined by \eqref{c3casea}, then takes the following form
\begin{equation}\eqlab{c3eqncasea}
\begin{aligned}
r_3 '&= -\frac{1}{1+\e^{\epsilon_3^{-1} v_3}+\e^{-\epsilon_3^{-1}(1+v_3)} }\left[-1+\e^{(\alpha-\beta)\epsilon^{-1}} \e^{\epsilon_3^{-1}(1+v_3)} +\e^{\epsilon_3^{-1} v_3} \right]r_3,\\
 v_3'&=\frac{\e^{-\epsilon_3^{-1}(1+v_3)}-\e^{(\m-\beta)\epsilon^{-1}} -\e^{-\epsilon_3^{-1}}}{1+\e^{\epsilon_3^{-1} v_3}+\e^{-\epsilon_3^{-1}(1+v_3)} }\\
 &+\frac{1}{1+\e^{\epsilon_3^{-1} v_3}+\e^{-\epsilon_3^{-1}(1+v_3)} }\left[-1+\e^{(\alpha-\beta)\epsilon^{-1}} \e^{\epsilon_3^{-1}(1+v_3)} +\e^{\epsilon_3^{-1} v_3} \right]v_3,\\
\epsilon_3 '&= \frac{1}{1+\e^{\epsilon_3^{-1} v_3}+\e^{-\epsilon_3^{-1}(1+v_3)} }\left[-1+\e^{(\alpha-\beta)\epsilon^{-1}} \e^{\epsilon_3^{-1}(1+v_3)} +\e^{\epsilon_3^{-1} v_3} \right]\epsilon_3.
\end{aligned}
 \end{equation}
Here we have as discussed above used a desingularization corresponding to division of the right hand side by $\epsilon_3$. Moreover, we have for convenience also multiplied the right hand side by
\begin{align*}
 \frac{\sum_{i=1}^6 e^{F_i(\beta-r_3,r_3 v_3)\epsilon^{-1}} }{1+\e^{\epsilon_3^{-1} v_3}+\e^{-\epsilon_3^{-1}(1+v_3)} },
\end{align*}
which is regular and nonzero on $D_3$. As above, we have left $\epsilon=r_3\epsilon_3$ as a constant where it makes sense. The system clearly extends smoothly to $\{r_3=0\}$ and $\{\epsilon_3=0\}$ provided that $v_3\ne 0$ and $v_3\ne -1$. In fact, these sets are by construction invariant. Within $\{r_3=0\}$ we have
\begin{equation}\nonumber
\begin{aligned}
 v_3'&=\frac{\e^{-\epsilon_3^{-1}(1+v_3)}- \e^{-\epsilon_3^{-1}}}{1+\e^{\epsilon_3^{-1} v_3}+\e^{-\epsilon_3^{-1}(1+v_3)} }\\
 &+\frac{1}{1+\e^{\epsilon_3^{-1} v_3}+\e^{-\epsilon_3^{-1}(1+v_3)} }\left[-1 +\e^{\epsilon_3^{-1} v_3} \right]v_3,\\
\epsilon_3 '&= \frac{1}{1+\e^{\epsilon_3^{-1} v_3}+\e^{-\epsilon_3^{-1}(1+v_3)} }\left[-1+\e^{\epsilon_3^{-1} v_3} \right]\epsilon_3.
\end{aligned}
 \end{equation}
 and therefore within the intersection $\{r_3=\epsilon_3=0\}$, we obtain
\begin{align*}
 v_3' &=\begin{cases}
         1 & \mbox{for $v_3<-1$},\\
         -v_3 & \mbox{for $-1<v_3<0$},\\
         v_3 & \mbox{for $v_3>0$},
        \end{cases}
\end{align*}
for any $v_3\ne -1$ and $v_3\ne 0$. In particular, $v_3'>0$ for these values of $v_3$.
%
The points
$v_3=0$ and $v_3=-1$ are degenerate (nonsmooth) and this stems from the lack of smoothness along $\mathcal E_{13}$ and $\mathcal E_{14}$, see \figref{pp}(a). We only resolve the degeneracy along $\mathcal E_{13}$, as $\mathcal E_{14}$ is not relevant for slow manifolds. For this we apply the blowup \eqref{bu1casea} to \eqref{c3eqncasea} (it turns out that multiplication by $\epsilon_3$ is not necessary in this case). We will study the dynamics in the chart \eqref{ci2casea} with $i=3$ in the following section.
\subsection{Analysis in the $\bar \epsilon_3=1$-chart}\seclab{bareps3}
Denote the vector-field associated with \eqref{c3eqncasea} by $V_3$. We then consider the pull-back of $V_3$ by the map $(r_3,v_{32},\rho_{32})\mapsto (r_3,v_3,\epsilon_3)$ defined by \eqref{ci2casea} with $i=3$. This gives the following equations
 \begin{equation}\eqlab{c32eqncasea}
\begin{aligned}
r_3 '&= -\left[-1+\e^{(\alpha-\beta)\epsilon^{-1}} \e^{\rho_{32}^{-1}}\e^{v_{32}} +\e^{v_{32}} \right] r_3,\\
 v_{32}'&=\rho_{32}^{-1}\e^{-\rho_{32}^{-1}} \left( \e^{v_{32}} -1-\e^{(\m-\beta)\epsilon^{-1}}\e^{\rho_{32}^{-1}}\right),\\
\rho_{32} '&= \left[-1+\e^{(\alpha-\beta)\epsilon^{-1}} \e^{\rho_{32}^{-1}} \e^{v_{32}} +\e^{v_{32}} \right]\rho_{32},
\end{aligned}
 \end{equation} after multiplying the right hand side by $1+\e^{v_{32}}+\e^{-\rho_{32}^{-1}} \e^{-v_{32}} $, which is regular and nonzero within a  fixed compact $D_{32}$ defined by $(r_3,v_{32},\rho_{32})\in [0,r_{30}]\times I\times [0,\rho_{320}]$, with $r_{30},\rho_{320}>0$ both sufficiently small and $I\subset \mathbb R$ a compact interval containing $0$. The system now extends smoothly to all of $D_{32}$. Within $\rho_{32}=0$, we have
 \begin{align*}
  r_3 ' &=-\left[\e^{v_{32}}-1\right]r_3,\\
  v_{32}'&=0.
 \end{align*}
It follows that any point $(0,v_{32},0)$ with $v_{32}\ne 0$ is partially hyperbolic. However, the line $v_{32}=\rho_{32}=0$, $r_3\ge 0$, is completely degenerate. Moreover, this degeneracy is clearly due to the flat factor $\e^{-\rho_3^{-1}}$ in the $v_{32}$-equation and consequently it cannot be resolved by the usual blowup method. We therefore use the method from \cite{kristiansen2017a} in the following section.

In conclusion, we have obtained a normally hyperbolic and attracting critical manifold on the blown up space, but only up until the blowup of $\mathcal E_{13}\times \{0\}$ (where the critical manifold loses hyperbolicity). We summarize the results of the section in \figref{ppabulocal2}.

 \begin{figure}
  \begin{center}
    {\includegraphics[width=0.475\textwidth]{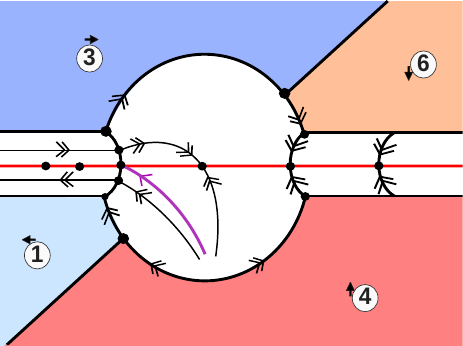}}
%
    \caption{Illustration of the results in the case (iii) on the blowup of $(u,v,\epsilon)=(\beta,0,0)$. The critical manifold is in red. It is normally hyperbolic and attracting everywhere except along the blowup of $\mathcal E_{13}$ (the left most cylinder). Notice that the third dimension coming out of the page corresponds to $\epsilon$. }\figlab{ppabulocal2}
        \end{center}
    \end{figure}

\section{Resolving the dynamics near $\mathcal E_{13}$ in the case (iii)}\seclab{flat}
To study \eqref{uv} near $\mathcal E_{13}:\,u\in[\alpha,\beta],v=0$, we consider $u\in [\alpha+\epsilon \xi,\beta-\epsilon \xi]$, for some fixed $\xi>0$, and set $v=\epsilon v_{2}$ with $v_{2}\in I$ where $I\subset \mathbb R$ is a fixed compact interval containing $0$. We will divide the analysis into (iii$1$) $\gamma<\alpha<\beta$ and (iii$2$) $\alpha<\gamma<\beta$.

In case (iii$1$), we define
\begin{align}\eqlab{qzdef}
 \begin{cases}
q = \e^{(u-\beta)\epsilon^{-1}},\\
z = \e^{(\alpha-u)\epsilon^{-1}}.
 \end{cases}
\end{align}
Notice that $q$ and $z$ are uniformly bounded as function of $u\in [\alpha+\epsilon \xi,\beta-\epsilon \xi]$. Also for $u=\beta+\mathcal O(\epsilon)$, we have $q=\mathcal O(1)$ and $z$ exponentially small, whereas for $u=\alpha+\mathcal O(\epsilon)$, we have $z=\mathcal O(1)$ and $q$ exponentially small. The idea is to augment $q$ and $z$ as variables, using
\begin{align}\nonumber
 \begin{cases}
\dot q = \epsilon^{-1} q\dot u,\\
\dot z = -\epsilon^{-1} z\dot u,
 \end{cases}
\end{align}
and define $v_{2}$ by $v=\epsilon v_{2}$, recall \eqref{vv32} and \eqref{c2casea}.
In this way, we can eliminate $u$ and obtain the following extended system from \eqref{uv}:
\begin{equation}\eqlab{qzv32}
\begin{aligned}
 \dot q &= q\left[v_{2}f(v_{2})+z \e^{v_{2}}\right],\\
 \dot z &= -z\left[v_{2}f(v_{2})+z \e^{v_{2}}\right],\\
 \dot v_{2} &=q \left(-\e^{-v_2}\left[v_{2}f(v_{2})+z \e^{v_{2}}\right]-z \e^{(\gamma-\alpha)\epsilon^{-1}} \right),
\end{aligned}
\end{equation}
after multiplication of the right hand side by $\epsilon \sum_{i=1}^6 \e^{ F_i(u,v)\epsilon^{-1}}$ so that
\begin{align}\eqlab{ueqnhere}
\dot u = \epsilon \left[v_{2}f(v_{2})+z \e^{v_{2}}\right].
\end{align}
Here we have introduced the smooth function $f:\mathbb R\rightarrow \mathbb R$ defined by
\begin{align}\eqlab{ffunc}
 f(x): = \begin{cases}
         1 & x=0,\\
         \frac{\e^{x}-1}{x} &x\ne 0.
        \end{cases}
\end{align}

We see that $(q,z,v_2)=(0,0,0)$ is fully nonhyperbolic for \eqref{qzv32}, since the linearization only has zero eigenvalues. We therefore introduce the following blowup transformation
\begin{align}\eqlab{casea2bu}
 \mu\ge 0,\, (\bar q,\bar z,\bar v_2)\in S^2\mapsto \begin{cases}
                                                     q = \mu \bar q,\\
                                                     z = \mu \bar z,\\
                                                     v_2 =\mu \bar v_2,
                                                    \end{cases}
\end{align}
and use a desingularization corresponding to division of the right hand side by $\mu$. We will demonstrate that this leads to a gain of hyperbolicity.
We consider the following charts
\begin{align}
 \bar q&=1:\quad \begin{cases}
                 q =\mu_1,\\
                 z =\mu_1 z_1,\\
                 v_2  =\mu_1 v_{21},
                \end{cases}\eqlab{c1a2}\\
                \bar z&=1:\quad \begin{cases}
                 q =\mu_2 q_2,\\
                 z =\mu_2,\\
                 v_2 =\mu_2 v_{22}.
                \end{cases}\eqlab{c2a2}
\end{align}
The change of coordinates are given by
\begin{align}\eqlab{ccfinal}
 \begin{cases}
  \mu_1 = \mu_2 q_2,\\
  z_1 = q_2^{-1},\\
  v_{21} = v_{22} q_2^{-1}.
 \end{cases}
\end{align}
\begin{remark}\remlab{invset}
By construction the equations \eqref{qzdef} define an invariant set for \eqref{qzv32}. On this set, the system \eqref{qzv32} is equivalent with \eqref{uv}. Notice in particular that on this invariant set we have $$z_1=\frac{z}{q}=\e^{(\alpha+\beta-2u)\epsilon^{-1}}.$$ Consequently, $z_1=\mathcal O(1)$ corresponds to $u=\frac{\alpha+\beta}{2}+\mathcal O(\epsilon)$ as $\epsilon\to 0$.
\end{remark}


In case (iii$2$), we will also augment
\begin{align}
 p=\e^{(\m-u)\epsilon^{-1}},\eqlab{pdef}
\end{align}
as a new variable
using $\dot p = - \epsilon^{-1} p \dot u$. Notice that $p$ is uniformly bounded for $u\in [\gamma+\epsilon \xi,\beta-\epsilon \xi]$, seeing that $\alpha<\gamma<\beta$ in the case (iii$2$). For $u\in [\alpha+\epsilon \xi,\gamma+\epsilon \xi]$, one could replace $p$ by $p^{-1}$, but it will suffice for our purposes to restrict attention to $$u\in [\tfrac{\gamma+\beta}{2}-\epsilon \xi,\beta-\epsilon \xi].$$ We therefore consider the following system:
\begin{equation}\eqlab{qzpv32}
\begin{aligned}
 \dot q &= q\left[v_{2}f(v_{2})+z \e^{v_{2}}\right],\\
 \dot z &= -z\left[v_{2}f(v_{2})+z \e^{v_{2}}\right],\\
 \dot p &=-p \left[v_{2}f(v_{2})+z \e^{v_{2}}\right],\\
 \dot v_{2} &=q \left(-v_2 f(-v_2)-p \right),
\end{aligned}
\end{equation}
using $\dot p = -\epsilon^{-1}p\dot u$.
Here $(q,z,p,v_2)=(0,0,0,0)$ is fully nonhyperbolic and we therefore use a blowup:
\begin{align}\eqlab{casea1bu}
 \mu\ge 0,\, (\bar q,\bar z,\bar p,\bar v_2)\in S^3\mapsto \begin{cases}
                                                     q = \mu \bar q,\\
                                                     z = \mu \bar z,\\
                                                     p = \mu \bar p,\\
                                                     v_2= \mu \bar v_2,
                                                    \end{cases}
\end{align}
and desingularization corresponding to division of the right hand side by $\mu$. We will demonstrate that this leads to a gain of hyperbolicity. For our purposes, it will suffice to consider the following chart
\begin{align}
 \bar q&=1:\quad \begin{cases}
                 q =\mu_1,\\
                 z =\mu_1 z_1,\\
                 p = \mu_1 p_1,\\
                 v_2  =\mu_1 v_{21}.
                \end{cases}\eqlab{c1a1}
\end{align}
Notice that by \eqref{qzdef} and \eqref{pdef}, we have
\begin{align}\eqlab{z1p1}
\begin{cases}
 z_1 = \frac{z}{q}=\e^{(\alpha+\beta-2u)\epsilon^{-1}},\\
 p_1 = \frac{p}{q}=\e^{(\m+\beta-2u)\epsilon^{-1}}.
 \end{cases}
 \end{align}
By assumption $\alpha<\gamma<\beta$ in the case (iii$2$), we conclude that
\begin{align}\eqlab{condz1p1}
z_1\ll p_1,
\end{align}
for $0<\epsilon\ll 1$ uniformly for $u\in [\tfrac{\gamma+\beta}{2}-\epsilon \xi,\beta-\epsilon \xi]$, $\xi>0$.
We will consider the two cases (iii$1$) and (iii$2$) successively in the following two sections.
\section{Case (iii$1$)}\seclab{casea2}
In this section, we consider \eqref{qzv32} with $\gamma<\alpha<\beta$ and apply the blowup \eqref{casea2bu}. We analyze the dynamics in the charts $\bar q=1$ and $\bar z=1$ successively.

\subsection{Chart $\bar q=1$}
In the chart $\bar q=1$, we consider the pull-back of the vector-field defined by \eqref{qzv32} by the mapping $(\mu_1,z_1,v_{21})\mapsto (q,z,v_2)$ given by \eqref{c1a2}. This gives
\begin{equation}\eqlab{c1a2eqns}
\begin{aligned}
 \dot \mu_1 &=\mu_1 \left[v_{21}f(\mu_1 v_{21})+z_1 \e^{\mu_1 v_{21}}\right],\\
 \dot z_1 &= -2 z_1 \left[v_{21}f(\mu_1 v_{21})+z_1 \e^{\mu_1 v_{21}}\right],\\
 \dot v_{21} &=-z_1 \e^{(\gamma-\alpha)\epsilon^{-1}} -\left(\e^{-\mu_1 v_{21}}+v_{21}\right)  \left[v_{21}f(\mu_1 v_{21})+z_1 \e^{\mu_1 v_{21}}\right],
\end{aligned}
\end{equation}
after desingularization, corresponding to division of the right hand side by $\mu_1$. We fix a compact set $D_{21}\subset \mathbb R^{3}$ consisting of $(\mu_1,z_1,v_{21})\in [0,\mu_{10}]\times [0,z_{10}]\times I$, with $\mu_{10}>0$, $z_{10}>0$, and $I\subset \mathbb R$ a fixed sufficiently large interval.

We see that for $z_1=0$ or $\e^{(\gamma-\alpha)\epsilon^{-1}}\to 0$, the set given by
\begin{align}\eqlab{m10eqn}
\left[v_{21}f(\mu_1 v_{21})+z_1 \e^{\mu_1 v_{21}}\right]  = 0,
\end{align}
defines a critical manifold, which is normally hyperbolic and attracting. By the implicit function theorem, the critical manifold can be written as a smooth graph over $(\mu_1,z_1)$:
\begin{align}\eqlab{critv21}
 v_{21} = z_1 m_{10}(\mu_1z_1),
\end{align}
with $m_{10}(0)= -1$. Here we have also used \eqref{ffunc}. Moreover, to obtain the form \eqref{critv21}, we have written \eqref{m10eqn} as
\begin{align*}
  z_1^{-1} v_{21}f((\mu_1z_1) (z_1^{-1} v_{21}))+ \e^{(\mu_1z_1) (z_1^{-1} v_{21})}=0,
\end{align*}
for $z_1\ne 0$.

Within the invariant set defined by $z_1=0$, we have
\begin{align*}
 \dot \mu_1 &=\mu_1 v_{21} f(\mu_1 v_{21}),\\
 \dot v_{21} &=-\left(\e^{-\mu_1 v_{21}}+v_{21}\right) v_{21}f(\mu_1 v_{21}).
 \end{align*}
Here the critical manifold $v_{21}=0$, $\mu_1\ge 0$, corresponds to the line $(r_3,v_{32},\rho_{32})=(r_3,0,0)$, $r_3\ge 0$, in the $\bar \epsilon_3=1$-chart of the system \eqref{c32eqncasea}, but now the line is normally hyperbolic, also at $\mu_1=0$. Moreover, there is also a hyperbolic saddle at $(\mu_1,v_{21})=(0,-1)$, with the linearization being
\begin{align*}
 \begin{pmatrix}
  -1 & 0 \\
  \star & 1
 \end{pmatrix}.
\end{align*}

On the other hand, within the invariant set $\mu_1=0$, we have
\begin{equation}\eqlab{z1v21}
\begin{aligned}
 \dot z_1 &= -2 z_1 \left[v_{21}+z_1\right],\\
 \dot v_{21} &=-z_1 \e^{(\gamma-\alpha)\epsilon^{-1}} -\left(1+v_{21}\right)  \left[v_{21}+z_1\right].
 \end{aligned}
 \end{equation}
Here we have a partially hyperbolic point at $(z_1,v_{21})=(0,0)$ and a hyperbolic source at $(z_1,v_{21})=(0,-1)$.
\begin{lemma}\lemmalab{q1a2lemma}
 The system \eqref{c1a2eqns}  has a two-dimensional slow manifold within $D_{21}$ for all $0< \epsilon\ll 1$ of the following graph form
 \begin{align*}
  v_{21} = z_1 \left( m_{10}(\mu_1z_1)+\e^{(\gamma-\alpha)\epsilon^{-1}} m_{11}(\mu_1,z_1,\e^{(\gamma-\alpha)\epsilon^{-1}})\right),
 \end{align*}
with $m_{11}$ smooth and satisfying
\begin{align*}
 m_{11}(0,z_1,0)=-\frac{1}{1+z_1}.
\end{align*}

\end{lemma}
\begin{proof}
The result follows from a simple calculation.

\end{proof}
A simple calculation shows that the reduced problem on the slow manifold is given by
\begin{equation}\nonumber
\begin{aligned}
 \dot \mu_1 &=\mu_1z_1\e^{(\gamma-\alpha)\epsilon^{-1}} \left[-\frac{1}{1+z_1}+\mathcal O(\mu_1,\e^{(\gamma-\alpha)\epsilon^{-1}})\right],\\
 \dot z_1 &= -2 z_1^2 \e^{(\gamma-\alpha)\epsilon^{-1}} \left[-\frac{1}{1+z_1}+\mathcal O(\mu_1,\e^{(\gamma-\alpha)\epsilon^{-1}})\right].
 \end{aligned}
\end{equation}

\begin{remark}
Notice also that from \lemmaref{q1a2lemma}, upon blowing down and restricting to the invariant set defined by \eqref{qzdef}, we obtain a slow manifold of \eqref{uv} of the graph form
\begin{align*}
 v &= \epsilon \mu_1 z_1 \left( m_{10}(\mu_1z_1)+\e^{(\gamma-\alpha)\epsilon^{-1}} m_{11}(\mu_1,z_1,\e^{(\gamma-\alpha)\epsilon^{-1}})\right)\\
 &=\epsilon \e^{(\alpha-u)\epsilon^{-1}} \left( m_{10}(\e^{(\alpha -u)\epsilon^{-1}})+\e^{(\gamma-\alpha)\epsilon^{-1}} m_{11}(\e^{(u-\beta)\epsilon^{-1}},\e^{(\alpha+\beta-2u)\epsilon^{-1}},\e^{(\gamma-\alpha)\epsilon^{-1}})\right),
\end{align*}
for $u\in [\frac{\alpha+\beta}{2}-\epsilon \xi,\beta-\epsilon \xi]$ with $\xi>0$ small. Here we have also used \eqref{c1a2} and \eqref{z1p1}.
\end{remark}
\subsection{Chart $\bar z=1$}
In the chart $\bar z=1$, we consider the pull-back of the vector-field defined by \eqref{qzv32} by the mapping $(\mu_2,q_2,v_{22})\mapsto (q,z,v_2)$ given by \eqref{c2a2}.  This gives
\begin{equation}\eqlab{c2a2eqns}
\begin{aligned}
\dot \mu_2 &= -\mu_2 \left[v_{22} f(\mu_2 v_{22})+\e^{\mu_2 v_{22}}\right],\\
\dot q_2 &= 2q_2 \left[v_{22} f(\mu_2 v_{22})+\e^{\mu_2 v_{22}}\right],\\
\dot v_{22} &=-q_2 \e^{(\gamma-\alpha)\epsilon^{-1}} +\left(-q_2 \e^{-\mu_2v_{22}}+v_{22} \right)\left[v_{22} f(\mu_2 v_{22})+\e^{\mu_2 v_{22}}\right],
\end{aligned}
\end{equation}
after desingularization, corresponding to division of the right hand side by $\mu_2$. We fix a compact set $D_{22}\subset \mathbb R^3$ consisting of $(\mu_2,q_2,v_{22})\in [0,\mu_{20}]\times [0,q_{20}]\times I$, with $\mu_{20}>0$, $q_{20}>0$, and $I\subset \mathbb R$ a fixed sufficiently large interval. We can take $q_{20}>0$ large enough so that $D_{22}$ when written in the $(\mu_1,z_1,v_{21})$-coordinates (using the change of coordinates, which are well-defined for $q_2>0$) intersects $D_{21}$ with nonzero interior.

We see that $(\mu_2,q_2,v_{22})=(0,0,0)$ is fully hyperbolic, with linearization
\begin{align*}
 \begin{pmatrix}
  -1 & 0 & 0\\
  0 & 2 & 0\\
  0 & \star & 1
 \end{pmatrix}.
\end{align*}
Moreover, for $q_2=0$ or $\e^{(\gamma-\alpha)\epsilon^{-1}}\to 0$ we have that the set given by
\begin{align*}
\left[v_{22} f(\mu_2 v_{22})+\e^{\mu_2 v_{22}}\right]=0,
\end{align*}
defines a normally hyperbolic and attracting critical manifold. By the implicit function theorem, the manifold is a graph over $\mu_2$ of the following form
\begin{align}\eqlab{critv22}
 v_{22} = m_{10}(\mu_2),\quad m_{10}(0)=-1,
\end{align}
for all $\mu_2\in [0,\mu_{20}]$, with $\mu_{20}>0$ small enough, recall \eqref{critv21}. The critical manifold corresponds to \eqref{critv21} upon change of coordinates.

\begin{lemma}\lemmalab{z1a2lemma}
 The system \eqref{c2a2eqns}  has a two-dimensional slow manifold within $D_{22}$ for all $0< \epsilon\ll 1$ of the following graph form
 \begin{align*}
  v_{22} = m_{10}(\mu_2)+q_2 \e^{(\gamma-\alpha)\epsilon^{-1}} m_{21}(\mu_2,q_2,\e^{(\gamma-\alpha)\epsilon^{-1}}),
 \end{align*}
with $m_{21}$ smooth and satisfying
\begin{align*}
 m_{21}(0,q_2,0)=-\frac{1}{1+q_2}.
\end{align*}
\end{lemma}
\begin{proof}
The result follows from a simple calculation.
\end{proof}
A simple calculation, shows that the reduced problem on the slow manifold is given by
\begin{equation}\nonumber
\begin{aligned}
 \dot \mu_2 &=-\mu_2q_2\e^{(\gamma-\alpha)\epsilon^{-1}} \left[-\frac{1}{1+q_2}+\mathcal O(\mu_2,\e^{(\gamma-\alpha)\epsilon^{-1}})\right],\\
 \dot q_2 &= 2 q_2^2 \e^{(\gamma-\alpha)\epsilon^{-1}} \left[-\frac{1}{1+q_2}+\mathcal O(\mu_2,\e^{(\gamma-\alpha)\epsilon^{-1}})\right].
 \end{aligned}
\end{equation}

\begin{remark}
Notice also that from \lemmaref{z1a2lemma}, upon blowing down and restricting to the invariant set defined by \eqref{qzdef}, we obtain a slow manifold of \eqref{uv} of the graph form
\begin{align*}
 v &= \epsilon \mu_2 \left( m_{10}(\mu_2)+q_2 \e^{(\gamma-\alpha)\epsilon^{-1}} m_{21}(\mu_2,q_2,\e^{(\gamma-\alpha)\epsilon^{-1}})\right)\\
 &=\epsilon \e^{(\alpha-u)\epsilon^{-1}} \left( m_{10}(\e^{(\alpha-u)\epsilon^{-1}})+\e^{(2u-\alpha-\beta)\epsilon^{-1}} \e^{(\gamma-\alpha)\epsilon^{-1}} m_{21}(\e^{(\alpha-u)\epsilon^{-1}},\e^{(2u-\alpha-\beta)\epsilon^{-1}},\e^{(\gamma-\alpha)\epsilon^{-1}})\right),
\end{align*}
for $u\in [\alpha+\epsilon \xi,\frac{\alpha+\beta}{2}+\epsilon \xi]$ with $\xi>0$ fixed and all $0<\epsilon\ll 1$. Here we have also used that $q_2=z_1^{-1}=\e^{(2u-\alpha-\beta)\epsilon^{-1}}$.
\end{remark}
We summarize the results of the section in \figref{casea2}.

 \begin{figure}
  \begin{center}
    {\includegraphics[width=0.475\textwidth]{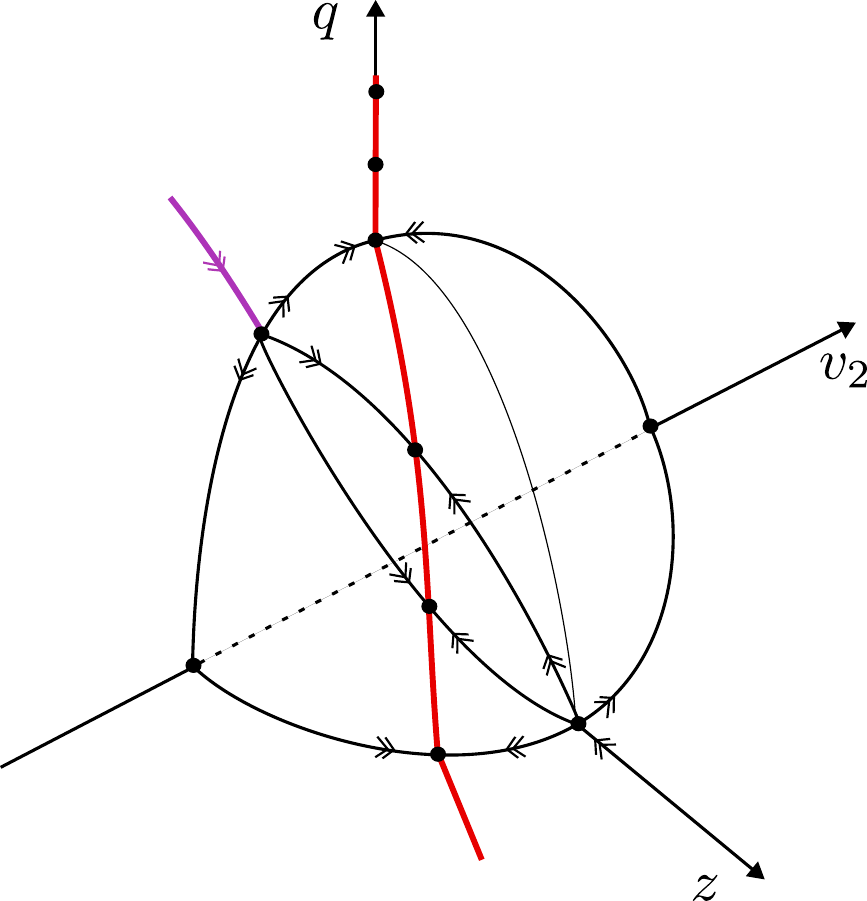}}
%
    \caption{Illustration of the results in the case (iii$1$) on the blowup of $(q,z,v_2)=(0,0,0)$ by \eqref{casea2bu}. The critical manifold for $\epsilon=0$ is in red. It is normally hyperbolic and attracting everywhere on the blowup sphere. The purple orbit within $\{z=0\}$ is a stable manifold of a hyperbolic saddle (at $(\mu_1,z_1,v_{21})=(0,0,-1)$ in the coordinates of chart $\bar q=1$) and it corresponds to the purple orbit in \figref{ppabulocal2} with the $\omega$-limit: $(r_3,v_{32},\rho_{32})=(0,0,0)$ in the chart $\bar \epsilon_3=1$, see \secref{bareps3}, upon restricting to the invariant set defined by \eqref{qzdef}. }\figlab{casea2}
        \end{center}
    \end{figure}

\section{Case (iii$2$)}\seclab{casea1}
In this section, we consider the pull-back of \eqref{qzpv32} by the mapping $(q,z,p,v_2)\mapsto (\mu_1,z_1,p_1,v_{21})$ defined by \eqref{c1a1}. This gives
\begin{equation}\eqlab{qzv32a1}
\begin{aligned}
 \dot \mu_1 &= \mu_1 \left[v_{21}f(\mu_1 v_{21})+z_1 \e^{\mu_1 v_{21}}\right],\\
 \dot z_1 &= -2z_1 \left[v_{21}f(\mu_1 v_{21})+z_1 \e^{\mu_1 v_{21}}\right],\\
 \dot p_1 &=-2p_1 \left[v_{21}f(\mu_1 v_{21})+z_1 \e^{\mu_1 v_{21}}\right],\\
 \dot v_{21} &= -v_{21} f(-\mu_1 v_{21})-p_1-\left[v_{21}f(\mu_1 v_{21})+z_1 \e^{\mu_1 v_{21}}\right]v_{21},
\end{aligned}
\end{equation}
after dividing the right hand side by $\mu_1$. Here $(\mu_1,z_1,p_1,v_{21})=(0,0,0,0)$ is partially hyperbolic for \eqref{qzv32a1} with eigenvalues $-1,0,0,0$ of the linearization.

Now, within $\mu_1=z_1=0$ we have
\begin{equation}\eqlab{p1v21eqn}
\begin{aligned}
 \dot p_1 &=-2p_1v_{21},\\
 \dot v_{21} &= -v_{21} -p_1-v_{21}^2,
\end{aligned}
\end{equation}
Here $(p_1,v_{21})=(0,0)$ is partially hyperbolic whereas $(p_1,v_{21})=(0,-1)$ is a hyperbolic unstable node. The flow of \eqref{p1v21eqn} is regular for $p_1>0$.

\begin{lemma}\lemmalab{center}
 There exists a three-dimensional center manifold of $(\mu_1,z_1,p_1,v_{21})=(0,0,0,0)$ of the following graph form
 \begin{align*}
  v_{21} =-p_1(1+\mathcal O(\mu_1,z_1,p_1)).
 \end{align*}
The submanifold within $\mu_1=z_1=0$:
\begin{align}\eqlab{cma1mu2z1eq0}
 v_{21} = -p_1(1+\mathcal O(p_1)),
\end{align}
is unique.
\end{lemma}
\begin{proof}
 Simple calculation. The fact that the submanifold within $\{\mu_1=z_1=0\}$ is unique follows from the fact that $(p_1,v_{21})=(0,0)$ is a nonhyperbolic saddle with the flow on the center manifold given by
 \begin{align*}
  \dot p_1 &=2p_1^2 \left(1+\mathcal O(p_1)\right).
 \end{align*}

\end{proof}
On the center manifold we find that
\begin{align*}
  \dot \mu_1 &= \mu_1 \left[-p_1(1+\mathcal O(\mu_1,z_1,p_1))+z_1(1+\mathcal O(\mu_1,z_1,p_1))\right],\\
 \dot z_1 &= -2z_1 \left[-p_1(1+\mathcal O(\mu_1,z_1,p_1))+z_1(1+\mathcal O(\mu_1,z_1,p_1))\right],\\
 \dot p_1 &=-2p_1 \left[-p_1(1+\mathcal O(\mu_1,z_1,p_1))+z_1(1+\mathcal O(\mu_1,z_1,p_1))\right].
\end{align*}
By \eqref{condz1p1}, we see that $p_1$ and $z_1$ are increasing whereas $\mu_1$ is decreasing.

The center manifold provides an extension of the slow manifold in the $\bar \epsilon=1$-chart, recall \secref{bareps1} and \lemmaref{bareps1lemma}, into $p_1=\mathcal O(1)$. Notice that $p_1=\mathcal O(1)$ corresponds to
\begin{align}\eqlab{u2def}
 u = \frac{\beta+\gamma}{2}+\epsilon \log x,\quad x=\mathcal O(1),\,x>0,
\end{align}
 cf. \eqref{z1p1} and \eqref{meps}.
Moreover, since $v=\epsilon v_{2}=\epsilon \mu_1 v_{21}$, we also obtain from \eqref{u2def} that
\begin{align}\eqlab{v21def}
 v = \epsilon \e^{\frac{\gamma-\beta}{2\epsilon}} y,\quad y:=xv_{21},
\end{align}
cf. \eqref{qzdef}. Finally, we emphasize that \eqref{u2def} and \eqref{v21def} agree with \eqref{u2v21def}.

In fact, if we view \eqref{u2def} and \eqref{v21def} as a change of coordinates $(x, y)\mapsto (u,v)$ for \eqref{uv}, then we obtain the following system:
\begin{align*}
 \dot x &=\e^{\frac{\gamma-\beta}{2\epsilon}}\left(x
 y+\mathcal O(\e^{(\alpha-\gamma)\epsilon^{-1}},\e^{\frac{\gamma-\beta}{2\epsilon}})\right)\\
 \dot y &=\e^{\frac{\gamma-\beta}{2\epsilon}}\left(-x y-1+\mathcal O(\e^{(\alpha-\gamma)\epsilon^{-1}},\e^{\frac{\gamma-\beta}{2\epsilon}})\right),
\end{align*}
after a regular change of time. Therefore by dividing the right hand side by $\e^{\frac{\gamma-\beta}{2\epsilon}}$ and letting $\epsilon\to 0$ we obtain \eqref{xyeqn}, repeated here for convenience:
\begin{equation}\eqlab{u2tildev21}
\begin{aligned}
 x'
  &=xy,\\
y' &=-xy-1.
\end{aligned}
\end{equation}
Here $x>0$ and $y\in \mathbb R$.
This system is equivalent with \eqref{p1v21eqn} for $p_1>0$ by the change of coordinates defined by $p_1=x^{-2}$, $v_{21}=x^{-1}y$ and a regular change of time. This follows from a simple calculation. We show the phase portrait of \eqref{u2tildev21} in \figref{phaseportrait}. The red orbit indicates the unique center manifold \eqref{cma1mu2z1eq0} (globalized by the flow), which in the $(x,y)$-coordinates take the following form:
\begin{align}\eqlab{cma1xy}
  y = -x^{-1}(1+\mathcal O(x^{-1})).
\end{align}
  \begin{figure}
  \begin{center}
  {\includegraphics[width=0.85\textwidth]{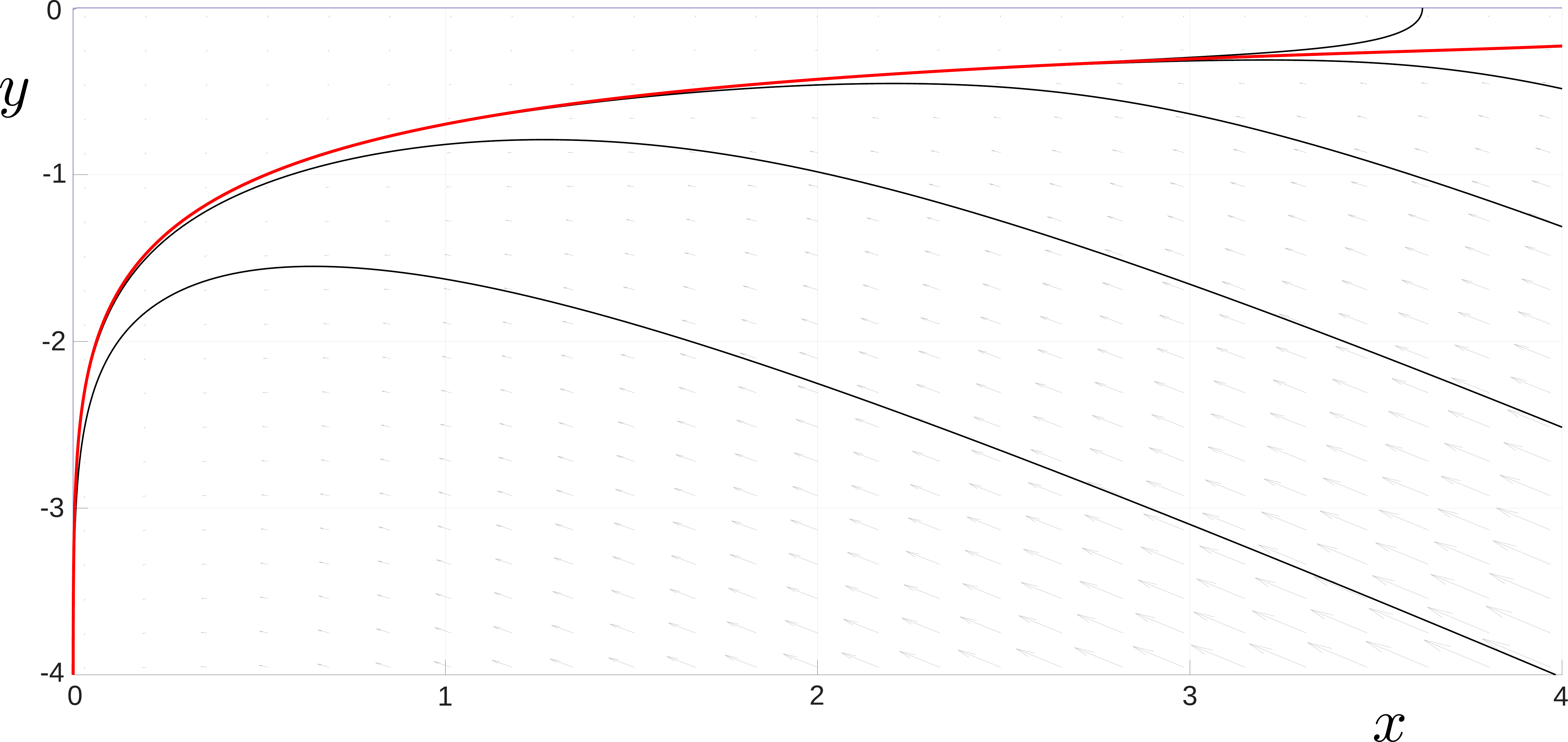}}
    \caption{Phase portrait  of \eqref{u2tildev21}. The orbit in red is the unique center manifold in \eqref{cma1xy}, globalized by the flow.}\figlab{phaseportrait}
        \end{center}
    \end{figure}

We summarize the results of the section in \figref{casea1}.

 \begin{figure}
  \begin{center}
    {\includegraphics[width=0.475\textwidth]{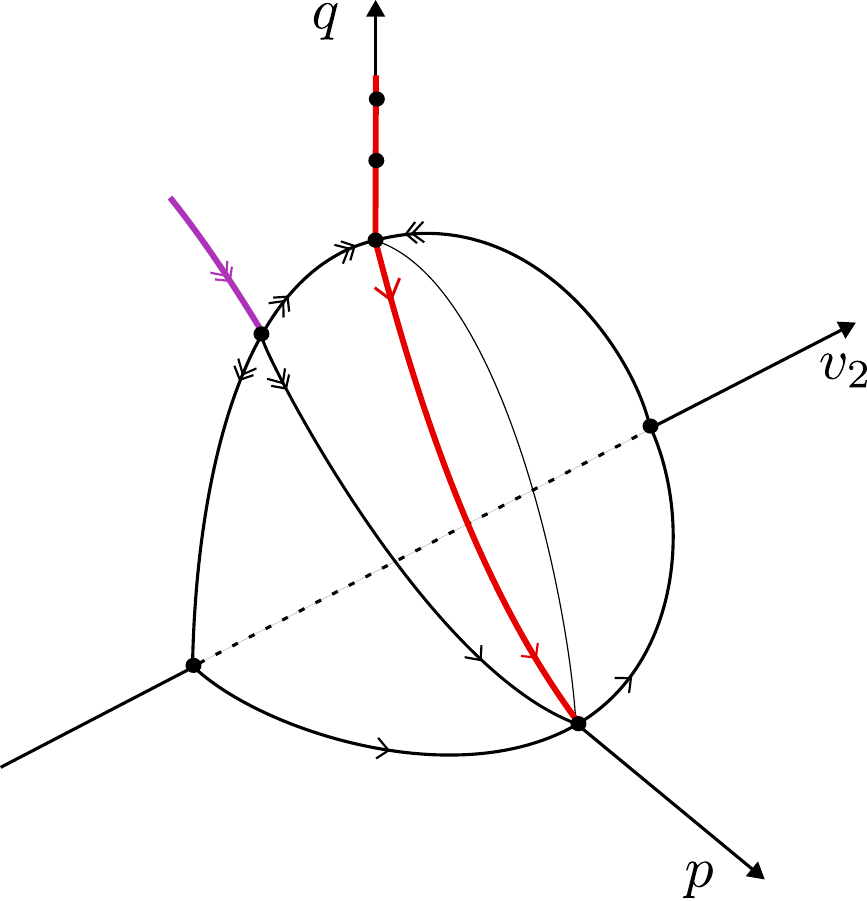}}
%
    \caption{Illustration of the results in the case (iii$2$) on the blowup of $(q,z,p,v_2)=(0,0,0,0)$ by \eqref{casea1bu}. Here we only illustrate $z=0$. On the sphere the flow is regular. In particular, the red orbit is the unique center manifold within $\{\mu_1=0\}$ of the partially hyperbolic point at $(\mu_1,p_1,v_{21})=(0,0,0)$, see \lemmaref{center}.}\figlab{casea1}
        \end{center}
    \end{figure}

%
\section{Completing the proof of \thmref{main}}\seclab{completing}
Consider first case (i). Then the results in \secref{casec} shows that there is a sequence of blowup and desingularizations (illustrated in \figref{ppcbulocal}) so that there is a globally defined critical manifold that is normally hyperbolic and attracting in compact and overlapping charts. Notice in particular by \remref{caseceps1extrem} that the unique invariant manifold from infinity is truly a slow manifold. It then follows from Fenichel's theory \cite{fen3} that the unique invariant manifold from infinity is a global slow manifold, recall also \remref{caseceps3extrem}.

Next, we consider case (ii). Instead of the blowup \eqref{bucasec1} used for (i), we now blowup $(u,v,\epsilon)=(\alpha,0,0)$ using the following spherical blowup transformation:
\begin{align}
 r\ge 0,\,(\bar u,\bar v,\bar \epsilon)\in S^2\mapsto \begin{cases}
                                                       u =\alpha+r\bar u,\\
                                                       v =r\bar v,\\
                                                       \epsilon =r\bar \epsilon.
                                                      \end{cases}
                                                      \end{align}
We consider the charts $\bar \epsilon=1$  and $\bar u=\pm 1$. In the latter charts, we blowup the sets corresponding to $\mathcal E_{46}\times \{0\}$ resp. $\mathcal E_{45}\times \{0\}$ as in \eqref{bu1casea} and \eqref{bu3casea}. The results are similar to the results in the case (i): We find a globally defined critical manifold that is normally hyperbolic and attracting in compact and overlapping charts. The properties in \remref{caseceps1extrem} and \remref{caseceps3extrem} also carry over and the statement in \thmref{main} on the case (ii) then follows from Fenichel's theory \cite{fen3}. The details of the slow flow are similar and left out for simplicity.


Finally, in the case (iii): From the analysis in the $\bar \epsilon=1$- and $\bar \epsilon_1=1$-charts, see \secref{bareps1} resp. \secref{bareps11}, and by \remref{caseaeps1extrem}, we conclude using Fenichel's theory that the unique invariant manifold from infinity is a slow manifold that extends, at least up until $\{u_2=-\xi\}$ in the chart $\bar \epsilon = 1$. From here we separate the cases (iii$1$) and (iii$2$). We start with (iii$1$). For this we will use the results of \secref{casea2}. 
%
First, we consider the $\bar \epsilon_3=1$-chart and the equations \eqref{c32eqncasea}. Based upon \eqref{c1a2} and \eqref{qzdef}, we apply the following blowup
\begin{align*}
 v_2 = \e^{(u-\beta)\epsilon^{-1}}v_{21} = \e^{-\rho_{32}^{-1}}v_{21},
\end{align*}
using $u=\beta-r_3$ in the second equality,
and a desingularization corresponding division of the right hand side by $\e^{(u-\beta)\epsilon^{-1}}=\e^{-\rho_{32}^{-1}}$. This gives a smooth, locally defined system in the  $(r_3,v_{21},\rho_{32})$-coordinates with $v_{21}=0$ being a normally hyperbolic and attracting critical manifold within both $\{\rho_{32}=0\}$ and $\{r_3=0\}$. We leave out the simple calculations. Then for $\rho_{32}\ge \rho_{320}>0$ we eliminate $\rho_{32}$ through the conservation of $\epsilon=r_3\rho_{32}$, writing $\rho_{32}=r_3^{-1} \epsilon$ (as in \remref{caseceps3extrem}) and $r_3=\beta -u$. This gives the following  system
\begin{equation}\eqlab{uv21eps}
\begin{aligned}
 \dot u &= \epsilon \left[v_{21}f(\mu_1 v_{21} )+z_1\e^{\mu_1v_{21}} \right],\\
 \dot v_{21} &=-z_1 \e^{(\gamma-\alpha)\epsilon^{-1}}-(\e^{-\mu_1 v_{21}}+v_{21})\left[v_{21}f(\mu_1 v_{21} )+z_1\e^{\mu_1v_{21}} \right],\\
 \dot \epsilon &=0,
\end{aligned}
\end{equation}
with $\mu_1=\e^{(u-\beta)\epsilon^{-1}}$ and $z_1=\e^{(\alpha+\beta-2u)\epsilon^{-1}}$.
This follows directly from \eqref{c1a2eqns} using \eqref{qzdef} and \eqref{c1a2}, recall also \remref{invset}. The system \eqref{uv21eps} has a smooth extension to $\epsilon=0$ for $u\in (\frac{\alpha+\beta}{2},\beta)$ by assumption $\gamma<\alpha<\beta$. Within $\{\epsilon=0,\,u\in (\frac{\alpha+\beta}{2},\beta)\}$, the equation $v_{21}=0$ defines a normally hyperbolic and attracting critical manifold for \eqref{uv21eps}. In order to resolve the nonsmoothness at $u=\frac{\alpha+\beta}{2}$ for $\epsilon=0$, we then apply the blowup transformation
\begin{align*}
 r\ge 0,\,(\bar u,\bar \epsilon)\in S^1\mapsto \begin{cases}
                                                  u = \frac{\alpha+\beta}{2}+r\bar u,\\
                                                  \epsilon = r \bar \epsilon.
                                                 \end{cases}
\end{align*}
Here we use the same symbols: $r\ge 0$ and $(\bar u,\bar \epsilon)$, as above for simplicity.
We will not need to perform desingularization. By working in the associated $\bar u=1$-chart, we can extend the critical manifold defined by $v_{21}=0$ as a normally hyperbolic and attracting set onto $r=0$. We leave out the simple calculations. Then in the $\bar \epsilon=1$-chart, setting
\begin{align}\eqlab{u2final}
u=\frac{\alpha+\beta}{2} + \epsilon u_2,
\end{align}
we perform a diffeomorphic change of coordinates $(u_2,v_{21})\mapsto (u_2,v_{22})$ defined by
\begin{align*}
 v_{22} = v_{21}\e^{2u_2}.
\end{align*}
Here we have used \eqref{ccfinal} and the fact that $z_1^{-1} =\e^{2u_2}$ when \eqref{u2final} holds.
In this way, we can continue the critical manifold $v_{22}=m_{10}(\mu_2)$, recall \eqref{critv22}, which equals $v_{22}=-1$ for $\mu_2=0$,  as a normally hyperbolic and attracting set to $u-\frac{\alpha+\beta}{2}<0$ and $\mathcal O(1)$ by working in the associated $\bar u=-1$-chart:
\begin{align*}
 \begin{cases}
u = \frac{\alpha+\beta}{2}-r_3,\\
\epsilon =r_3\epsilon_3.
\end{cases}
\end{align*}
For $r_3\ge r_{30}>0$ with $r_{30}>0$ fixed small enough, we eliminate $\epsilon_3=r_3^{-1}\epsilon$. In this way, we obtain the following system
\begin{align*}
 \dot u &=\epsilon \left[v_{22}f(\mu_2 v_{22})+\e^{\mu_2v_{22}}\right]\\
 \dot v_{22}&=-q_2 \e^{(\gamma-\alpha)\epsilon^{-1}}+(-q_2 \e^{-\mu_2 v_{22}} +v_{22})\left[v_{22}f(\mu_2 v_{22})+\e^{\mu_2v_{22}}\right],\\
 \dot \epsilon &=0.
\end{align*}
with $\mu_2=\e^{(\alpha-u)\epsilon^{-1}}$ and $q_2=\e^{(2u-\alpha-\beta)\epsilon^{-1}}$,
after desingularization.
This follows from \eqref{c2a2eqns}. This system has a smooth extension to $\epsilon=0$ for $u\in (\alpha,\frac{\alpha+\beta}{2})$ and on this set, we see that $v_{22}=-1$ defines a normally hyperbolic and attracting critical manifold. (Indeed for $\mu_2=q_2=0$, the $v_{22}$-equation reduces to $\dot v_{22}=v_{22}(1+v_{22})$, so that $v_{22}=-1$ is attracting whereas $v_{22}=0$ is repelling.) To proceed further, we blowup $(u,v,\epsilon)=(\alpha,0,0)$ as follows
\begin{align}\eqlab{bufinal}
 r\ge 0,\,(\bar u,\bar v,\bar \epsilon)\in S^2\mapsto \begin{cases}
                                                       u = \alpha+r\bar u,\\
                                                       v = r\bar v,\\
                                                       \epsilon =r\bar \epsilon.
                                                      \end{cases}
\end{align}
The analysis is similar to the analysis in the case (i) of the blowup of $(u,v,\epsilon)=(\gamma,0,0)$. In the  associated $\bar u=\pm 1$-charts, we blowup the resulting degeneracies along $v_1=0$, $\epsilon_1=0$ resp. $v_3=-1$, $\epsilon_3=0$ (stemming from the degeneracies along $\mathcal E_{13}$ resp. $\mathcal E_{12}$, recall \figref{pp}(c))  through cylindrical blowups, see \eqref{bu1} and \eqref{bu3}. In the associated $\bar \epsilon_1=1$-chart, see \eqref{ci2} with $i=1$, we then use the following blowup transformation
\begin{align*}
 v_2 = \e^{(\alpha-u)\epsilon^{-1}} v_{22} = \e^{-\rho_{12}^{-1}}v_{22}.
\end{align*}
based upon \eqref{c2a2} and \eqref{qzdef}. Here we have used $u=\alpha+r_1$ in the second equality.
This allows us to extend the critical manifold $v_{22}=m_{10}(\mu_2)$, as a normally hyperbolic and attracting set into a compact set of the $\bar \epsilon=1$-chart associated with the blowup \eqref{bufinal}. The analysis from here is then identical to case (i) and is therefore left out. In this way, we have proven that in the case (iii$1$) there is a sequence of blowups and desingularizations of the extended system \eqref{uvext} so that the blown up vector-field has a globally defined normally hyperbolic and attracting critical manifold. Fenichel's theory (along with the uniformity properties, similar to \remref{caseceps3extrem} and \remref{caseaeps1extrem}, along the blowup of $\mathcal E_{12}\times \{0\}$ and $\mathcal E_{46}\times \{0\}$)  then completes the proof for case (iii$1$).

Now we finally turn to the case (iii$2$). It follows from the analysis of \secref{casea1} that the slow manifold from infinity extends into $x=\mathcal O(1)$ and $y=\mathcal O(1)$ for all $0<\epsilon\ll 1$ (through the center manifold in \lemmaref{center} restricted to the invariant set defined by \eqref{qzdef}, \eqref{pdef} and \eqref{c1a1}), recall \eqref{u2def} and \eqref{v21def}. Here the flow is regular for $0<\epsilon\ll 1$. Therefore it is not slow-fast. The regular orbit connects with the attracting slow manifold stemming from $\mathcal E_{12}$. This (along with repeating the arguments in the case (iii$1$) on the blowup of $(u,v,\epsilon)=(\alpha,0,0)$) proves the claim (on the nonexistence of global slow manifold) for (iii$2$).

\section{Discussion}\seclab{conclude}
In this section, we first discuss the nongeneric scalings (corresponding to the boundary of our domains in the $(\alpha,\beta,\gamma)$-space) that we have not addressed in our tropical analysis. This is particularly relevant since the classical reduction methods (QSSA, tQSSA, rQSSA) lie on these boundaries. Subsequently, we compare our results with related results in the extensive literature on reduction methods for Michaelis--Menten-type reactions.
{%
For the sake of this discussion, we consider \eqref{mm-dim} and write $$K_S= \frac{k_{-1}}{k_1}, \quad K=\frac{k_2}{k_1}, \quad K_M=K_S+K,$$ for the well-known dissociation, catalytic and Michaelis constants. We can then express our nondimensionalized parameters $(\tilde \alpha,\tilde \beta,\tilde \gamma)$ as follows
\begin{align}\eqlab{abct}
\tilde \alpha=\frac{K_S}{s_0}, \quad \tilde \beta=\frac{e_0}{s_0}, \quad \tilde \gamma=\frac{K}{s_0}.\end{align} The conditions (i), (ii), (iii$1$), (iii$2$), see the start of \secref{main}, only compare $\alpha,\beta,\gamma$, recall \eqref{abc}. Hence they only involve the ratios of the three rates $k_{-1}$, $k_2$ and $k_1e_0$, and {\textit{do not depend on $s_0$}}.

The generic cases (i), (ii), (iii$1$) and (iii$2$) are open cones in the space of exponents $(\alpha,\beta,\gamma)\in\mathbb R^3$. They are separated by four \emph{boundaries}, i.e.\ half-planes on which two exponents coincide:
\begin{enumerate}
\item[(B1)] $\alpha=\gamma>\beta$, between (i) and (ii);
\item[(B2)] $\alpha=\beta>\gamma$, between (ii) and (iii$1$);
\item[(B3)] $\alpha=\gamma<\beta$, between (iii$1$) and (iii$2$);
\item[(B4)] $\beta=\gamma>\alpha$, between (i) and (iii$2$).
\end{enumerate}
All four boundaries meet along the line $\alpha=\beta=\gamma$. The remaining two half-planes on which two exponents coincide, $\alpha=\beta<\gamma$ and $\beta=\gamma<\alpha$, {belong to the} cases (i) and (ii), respectively.
The classical regimes lie on the boundaries: the standard QSSA ($e_0\ll  K_S\sim K$) on (B1), the total QSSA (tQSSA) ($K\ll e_0\sim K_S$) on (B2) and the reverse QSSA (rQSSA) ($K_S\sim K\ll e_0$) on (B3). Finally, (B4) corresponds to $K_S\ll e_0\sim K$.

On a boundary (B1)-(B4) some of the tropical lines in \figref{pp} coalesce and  the exponential scaling \eqref{abc} only fixes the order of magnitude of the two coinciding parameters. The $\mathcal O(1)$ values of their ratios:
\begin{align*}
\frac{\tilde \gamma}{\tilde \alpha}=\frac{k_2}{k_{-1}}\ \text{on (B1), (B3)},\qquad \frac{\tilde \alpha}{\tilde \beta}=\frac{k_{-1}}{k_1e_0}\ \text{on (B2)},\qquad \frac{\tilde \gamma}{\tilde \beta}=\frac{k_2}{k_1e_0}\ \text{on (B4)}.
\end{align*}
can be viewed as unfolding parameters connecting the regimes (i), (ii), (iii1), (iii2) across the boundaries.
We are convinced that the methods of the paper can be used to study these remaining cases, but leave this to future work.

}

{%
\subsection{Comparison with previous work}

\textbf{Classical results}. Segel \cite{segel1988a} and Segel--Slemrod \cite{segel1989a} derive the small parameter $\varepsilon_{\mathrm{SS}}=e_0/(K_M+s_0)$ from time scale estimates and prove convergence to the Michaelis--Menten equation as $\varepsilon_{\mathrm{SS}}\to 0$. Their main concern is the initial transient from $(s_0,0)$, which we do not consider. In our variables $\varepsilon_{\mathrm{SS}}=\tilde \beta/(1+\tilde \alpha+\tilde \gamma)$, so $\varepsilon_{\mathrm{SS}}\ll1$ if and only if $\beta<\max(\alpha,\gamma,0)$. This condition combines our classification with the position $u=0$ of the initial value. Their convergence theorem, with $s_0/K_M$ and $k_{-1}/k_2$ bounded, covers the case (i) and the boundary (B1). Their local criterion that the ratio of the time scales $t_C/t_S=Ke_0/(K_M+s)^2$ be small, evaluated at the current substrate concentration \cite[Sec.~4.3]{segel1989a}, becomes
\begin{align*}
\frac{t_C}{t_S} = \frac{\tilde \beta\tilde \gamma}{(\tilde \alpha+\tilde \gamma+\tilde u)^2}\ll1,
\end{align*}
in our scalings. Here we have used \eqref{tueqn} and \eqref{abct}.
It fails for $\tilde \alpha +\tilde \gamma +\tilde u\lesssim\sqrt{\tilde \beta\tilde \gamma}$, which in case (iii$2$): $\tilde \beta\gg \tilde \gamma\gg \tilde \alpha$ precisely corresponds to the scale \eqref{u2v21def} where the manifold from infinity ceases to be slow. In their treatment of the rQSSA, Segel and Slemrod exclude $k_2\ll k_{-1}$ because the rQSSA and the QSSA curves then coincide. In our terms, this excludes the case (iii$1$) and leaves (iii$2$) and (B3).

\textbf{Discrete scalings and Tikhonov--Fenichel parameter values}.
Lapuz and Wechselberger \cite{lapuz} catalogue the reductions of \eqref{mm} obtained by assigning each parameter $(\tilde \alpha,\tilde \beta,\tilde \gamma)$ one of the orders $\mathcal O(\varepsilon)$, $\mathcal O(1)$, $\mathcal O(\varepsilon^{-1})$, and compute the reduced flows coordinate-independently by the projection method of \cite{wechselberger2020a}. With $\varepsilon=\e^{-\epsilon^{-1}}$ these configurations are the lattice points $\{-1,0,1\}^3$ in our space of exponents. Of the $3^3=27$ configurations, $12$ belong to the generic cases of \thmref{main} ($5$ to (i), $5$ to (ii), one each to (iii$1$) and (iii$2$)), and $15$ lie on boundaries. The variables are kept on the single scale $\tilde u,\tilde v=\mathcal O(1)$, i.e.\ near the tropical point $(u,v)=(0,0)$. Critical manifolds consisting of horizontal and vertical branches {(whose intersections are not normally hyperbolic)} correspond to tropical edges and vertices, which our blowups resolve; for example, the corner in \cite[case R.2a.i]{lapuz}, which belongs to the case (iii$2$), is the vertex $(\alpha,0)$. Loss of hyperbolicity is resolved in the thesis \cite{lapuzphd}.

\textcolor{black}{We emphasize that the reference \cite{lapuz} also considers the reversible case of the Michaelis-Menten, which has a fourth parameter $\widetilde \delta$, see \cite[Eq. (2.8)]{lapuz}. This leads to $67$ additional cases in their discrete scalings with $\tilde u,\tilde v=\mathcal O(1)$, see \cite[Fig. 5]{lapuz}. We believe that it is also possible to study the reversible case with our tropical methods (setting $\widetilde \delta= \e^{\delta\epsilon^{-1}}$), but leave this to future work.}

Walcher and coworkers characterize the small parameters through Tikhonov--Fenichel parameter values \cite{goeke2015a}, at which the system has a normally hyperbolic manifold of equilibria, and compute the corresponding reductions \cite{noethen2011a,goeke2014a}; see also \cite{lax2020a} for the relation to scalings. For the irreversible system these values are $e_0=0$, $k_1=0$ and $k_2=0$ (three different points) \cite{goeke2015a}, with the rQSSA as the degenerate subcase $k_{-1}=k_2=0$ \cite{eilertsen2020a}. In our space of exponents they are points at infinity, each approached along one direction:
\begin{itemize}
\item $e_0=0$ is $\beta\to-\infty$, and $k_1=0$ is $\alpha,\gamma\to\infty$ with $\alpha-\gamma$ fixed (up to a rescaling of time). Both lie at the end of (i), (B1) and (ii).
\item $k_2=0$ is $\gamma\to-\infty$, at the end of (ii), (B2) and (iii$1$).
\item $k_{-1}=k_2=0$ is $\alpha,\gamma\to-\infty$, at the end of (iii$1$), (B3) and (iii$2$), depending on the relative rates.
\end{itemize}
 A Tikhonov--Fenichel reduction is valid on compact subsets of the phase space, i.e.\ on one scale of $\tilde u$.

Eilertsen and Schnell \cite{eilertsen2020a} identify $\eta=e_0/K_M$ as the natural small parameter of the QSSA. Like our classification, it does not depend on $s_0$, and $\eta\ll1$ is the condition $\beta<\max(\alpha,\gamma)$. The anti-quasi-steady-state conditions of \cite{eilertsen2022b} ($k_1e_0\approx k_2$, $k_{-1}\ll k_2$, no time scale separation near the steady state) describe the boundary (B4). For the open reaction (where there is a constant inflow of substrate), Eilertsen, Roussel, Schnell and Walcher \cite{eilertsen2021a} single out a distinguished trajectory from infinity as a candidate for a global slow manifold; without inflow it is our manifold from infinity.
}

\textbf{Computational singular perturbation}.
Patsatzis and Goussis \cite{patsatzis2019a,patsatzis2023a} use the algorithmic tools of computational singular perturbation (CSP) to derive a reduced model and local criteria for the validity of the quasi-steady-state and partial-equilibrium approximations along trajectories, and test them on a set of parameter cases. Like us, they are concerned with uniform approximation of slow manifolds in the whole phase space.  Their reduced model is based on the curve $c=e_0s/(K_S+FK+s)$, where $F\in(0,1)$ depends on $(s,c)$ through the eigenvalues of the Jacobian; it contains the standard QSSA ($F=1$), the reverse QSSA and the partial-equilibrium approximation as special or limiting cases. Their criteria are formulated in terms of the state-dependent quantities $\mu=K/(K_S+s)$ and $\nu=(e_0-c)/(K_S+s)$. At the steady state, $(\mu,\nu)=(k_2/k_{-1},k_1e_0/k_{-1})$, and the three regions of their $(\mu,\nu)$-plane separated by the lines $\nu=1+\mu$ and $\mu=1+\nu$ become our cases (iii), (ii) and (i); all their test cases fit into our classification. Along a trajectory, however, $(\mu,\nu)$ varies. On the scale \eqref{u2v21def}, with $(x,y)$ in a compact set, one has $\mu\approx\sqrt{K/e_0}\,x^{-1}\ll 1$ and $\nu\approx -y/x$ in the case (iii$2$), so this scale corresponds to $\mu\approx 0$, $\nu=\mathcal O(1)$, the left end of the line $\nu=1+\mu$, which is crossed at $x=-y$. The trajectories of \cite[case F]{patsatzis2019a}, which lies on (B3), and of \cite[case D]{patsatzis2023a}, which lies in (iii$2$), cross this line there, see \cite[Fig.~6]{patsatzis2019a} and \cite[Fig.~8]{patsatzis2023a}. The structural point that our analysis adds is that on this scale the dynamics is a regular perturbation of \eqref{xyeqn}, although the ratio of the eigenvalues of the Jacobian is small. 
A small ratio of eigenvalues therefore does not imply the existence of a slow manifold; cf.\ also \cite{eilertsen2022a}. In the case (iii$1$), the trajectory crosses $\nu=1+\mu$ at $s\approx\sqrt{e_0K_M}\approx\sqrt{e_0K_S}$, again with $\mu\ll 1$ and $\nu\approx 1$, but there the manifold from infinity is slow. Local criteria of this type therefore cannot distinguish (iii$1$) from (iii$2$).

\textbf{Tropical equilibration}.
Tropical methods for the reduction of chemical reaction networks were introduced by Radulescu and coworkers \cite{noel2014a,soliman2014a,desoeuvres2022a}. Parameters and concentrations are assigned orders of magnitude $\varepsilon^{a}$ with rational $a$ and fixed $\varepsilon$, and a concentration is called (partially) tropically equilibrated if, in its equation, the dominant positive and the dominant negative monomials have the same order. Equilibration is a necessary condition for a variable to be fast and relaxed. Together with a permanency assumption it leads to truncated reduced models. For the Michaelis--Menten reaction, \cite{soliman2014a} determine all branches of equilibrations in terms of the orders of $k_{-1}$, $k_2$ and $k_1e_0$, the same comparisons that define our cases. Our approach adds three things:
\begin{itemize}
\item the tropical limit is a dynamical (piecewise smooth) system \cite{kristiansen2026a}, with real exponents and $\epsilon\to0$;
\item the transitions between branches take place at tropical vertices, which we resolve by blowup;
\item existence and non-existence of slow manifolds are proved by GSPT.
\end{itemize}

By \secref{trop}, slow manifolds can only lie near the edges $\mathcal E_{45}$, $\mathcal E_{46}$ (cases (i), (ii)) and $\mathcal E_{12}$, $\mathcal E_{13}$, $\mathcal E_{46}$ (case (iii)). Each of these edges is a tropical equilibration in the sense of \cite{noel2014a,soliman2014a}: two monomials of opposite sign in the same equation balance. Along $\mathcal E_{12}$, $\mathcal E_{45}$ and $\mathcal E_{46}$ only one equation is equilibrated, which is a partial tropical equilibration in the sense of \cite{desoeuvres2022a}, see \figref{pp}. The edge $\mathcal E_{13}$ is different. Here $-\tilde u$ and $\tilde u\tilde v$ (terms $1$ and $3$, resp.) balance along $v=0$, with the corresponding flow vectors being tangent to the edge (tropical nullcline sliding, cf. \cite{kristiansen2026a}). At the same, the $v$-equation is also equilibrated (due to the balance of terms $4$ and $6$) on $v=0$ once $u>\max(\alpha,\gamma)$, see the red dotted curve in \figref{pp}(c). Tropically, (iii$1$) and (iii$2$) therefore look the same on $\mathcal E_{13}$: in both cases both equations are equilibrated there, which is the branch \cite[eq.~(23)]{soliman2014a}, for which a reduced model is proposed. Nevertheless, there is a slow manifold along $\mathcal E_{13}$ in the case (iii$1$), whereas in the case (iii$2$) the scale \eqref{u2v21def} on this edge carries regular, non-slow dynamics. Whether a slow manifold exists along a tangential edge is decided by subdominant (flat) terms, which become visible only after blowup (\secref{flat}). Tropical equilibration is therefore necessary but not sufficient for a slow manifold.

\section{Acknowledgement}
The authors thank TU Vienna (October 2025) and DTU (November 2026) for their hospitality during research visits. The first author was funded by Danish Research Council (DFF) grant 4283-00014B.

The authors used Anthrophic's Claude Opus 5.5 in the final stages of the research, i.e. for  proof-reading, for input on references, and for the comparison with previous work. In particular, Claude fixed some typos and wrote (under our strict guidance) a first draft of \secref{conclude}. In all cases, we checked the output and completed the final version. The authors take full responsibility for the content.
\bibliography{refs_tropicalization_mm}
\bibliographystyle{plain}
\end{document}